\documentclass[global,final]{svjour}[13.2.2006]
\journalname{European Journal of Combinatorics}

\usepackage {latexsym} \usepackage {amsmath} \usepackage {amsfonts} 
\usepackage {amssymb}
\usepackage {times}
\usepackage {epsf}
\usepackage {graphicx}
\usepackage {diagrams}
\usepackage {paralist}
\usepackage {ifpdf}
\smartqed

\ifpdf
\else
\fi

\let\phi\varphi
\def\?#1{\par{\bf ??? #1 !!!}\par}
\def\sss{\scriptscriptstyle}
\def\zet{{\mathbb Z}}

\def\id{\mathord{\mathrm {id}}}

\def\Ker{\mathop{\mathrm{Ker}}}

\DeclareMathAccent{\myarrow}{\mathord}{letters}{"7E}
\def\vecz{{\myarrow 0}}

\def\leh{\mathrel{\preccurlyeq_{h}}}
\def\lh{\mathrel{\prec_{h}}}

\def\leTT{\mathrel{\preccurlyeq}}
\def\geTT{\mathrel{\succcurlyeq}}
\def\lTT{\mathrel{\prec}}

\def\eqTT{\mathrel{\approx}}
\def\eqh{\mathrel{\approx_h}}
\newcommand{\chiTT}[1][]{\chi_{\sss TT_{#1}}}

\DeclareMathAccent{\myarrow}{\mathord}{letters}{"7E}
\def\orC{\overrightarrow C}
\def\orG{\overrightarrow G}
\def\orH{\overrightarrow H}
\def\orK{\overrightarrow K}
\def\barC{\overline C}
\def\symK{\overleftrightarrow{K}}  
\def\symG{\overleftrightarrow{G}}  
\def\symH{\overleftrightarrow{H}}  

\def\calG{{\cal G}}
\newcommand{\TT}[1][]{\mathrel{\xrightarrow{\ifx @#1@ TT\else TT_{#1} \fi}}}
\newcommand{\cc}{\mathrel{\xrightarrow{cc}}}
\newcommand{\nTT}[1][]{\thickspace\not\negthickspace\xrightarrow{\ifx @#1@ TT\else TT_{#1} \fi}}
\newcommand{\dunion}{\mathop{\dot\cup}}

\def\hom{\mathrel{\xrightarrow{hom}}}
\def\nhom{\mathrel{\thickspace\thickspace\not\negthickspace\negthickspace\xrightarrow{hom}}}
\def\floor #1{\left\lfloor #1 \right\rfloor}

\def\edge{\orK_2}
\def\sym{\mathop{\Delta}}

\begin{document}
\title{On tension-continuous mappings}
\author{
	Jaroslav Ne\v set\v ril
\and
	Robert \v S\'amal
}

\institute{Institute for Theoretical Computer Science (ITI)
   \thanks{Institute for Theoretical Computer Science is supported as project
    1M0021620808 by Ministry of Education of the Czech Republic.} 
 \\
 Charles University, Malostransk\'e n\'a\-m\v es\-t\'\i{}~25,
 118$\,$00 Prague, Czech Republic.
 \\
 \email {\{nesetril,samal\}@kam.mff.cuni.cz} }
\date{}
\maketitle
\begin{abstract}
Tension-continuous (shortly $TT$) mappings are mappings between 
the edge sets of graphs. They generalize graph homomorphisms.
From another perspective, tension-continuous mappings are dual to the notion
of flow-continuous mappings and the context of nowhere-zero flows
motivates several questions considered in this paper.

Extending our earlier research we define new constructions and
operations for graphs (such as graphs $\Delta_M(G)$) and give evidence
for the complex relationship of homomorphisms and $TT$ mappings.
Particularly, solving an open problem, we display pairs of
$TT$-comparable and homomorphism-incomparable graphs with
arbitrarily high connectivity.

We give a new (and more direct) proof of density of $TT$ order and
study graphs such that $TT$~mappings and homomorphisms from them coincide; 
we call such graphs homotens.
We show that most graphs are homotens, on the other hand 
every vertex of a nontrivial homotens graph is contained in a triangle. 
This provides a justification for our construction of homotens graphs.
\end{abstract}


\leftline{{\bfseries MSC\enspace}\enspace 05C15, 05C25, 05C38}

\section*{Contents}

\let\noi\noindent
\noi  1. Introduction 
  
\noi  2. Basic definitions

\noi  3. Examples

\noi  4. Left homotens

\noi  5. Right homotens

\noi  6. Density

\noi  7. Remarks

\section{Introduction}
\label{sec:intro}

It is a traditional mathematical theme to study the question when
a map between the sets of substructures is induced (as a lifting)
by a mapping of underlying structures. 
In a combinatorial setting (and as one of the simplest instances of
this general paradigm) this question takes the following form:

\begin{question}
\label{q:basic}
Given undirected graphs~$G$, $H$ and a mapping $g:E(G)\to E(H)$ does
there exist a mapping $f:V(G)\to V(H)$ such that 
$g(\{x,y\}) = \{f(x),f(y)\}$ for every edge $\{x,y\} \in E(G)$? 
\end{question}

In the positive case we say that \emph{$g$~is induced by~$f$}. 
It is easy to see that such mapping $f$ is a homomorphism $G \hom H$ and 
that to each homomomorphism corresponds exactly one induced mapping~$g$.
Thus Question~\ref{q:basic} asks which mappings~$g$ between edge sets
are induced by a homomorphism. 
Various instances of this problem were considered for example by Whitney
\cite{Whit1}, 
the first author \cite{Nes-derivative}, Kelmans \cite{Kelmans}, 
and by Linial, Meshulam, and Tarsi~\cite{LMT}. 
More recently, DeVos, Ne\v set\v ril, and Raspaud \cite{DNR} isolated
the following necessary condition for a mapping~$g : E(G) \to E(H)$
to be induced by a homomorphism. 

\begin{equation}\label{eq:def}
 \mbox{For every cut~$C \subseteq E(H)$ the set $g^{-1}(C)$ is a cut of~$G$.}
\end{equation}

Here, a \emph{cut} means the edge set of a spanning bipartite induced subgraph.
It is natural to call any mapping~$g$ satisfying condition~(\ref{eq:def})
a \emph{cut-continuous} mapping~$G \cc H$. Cut-continuous
mappings extend and generalize the notion of a homomorphism 
and the relationship of these two notions is the central theme of
this paper. We provide evidence in both directions. 
We present various examples of cut-continuous mapping that are not induced, 
in particular in Proposition~\ref{cnhomotens} we construct such mappings
between highly connected graphs, thereby solving a problem from our
previous paper~\cite{NS-TT1}. On the other hand, as described
in Section~\ref{sec:lh}, for most of the graphs all
cut-continuous mappings are induced.

Cut-continuous mappings were defined and investigated in~\cite{DNR,NS-TT1} 
in the  more general context of nowhere-zero flows and circuit
covers. As such, the tension-continuous mappings (being duals of
flow-continuous mappings) have deep combinatorial meaning. For
example, for a cubic graph~$G$ the number of cut-continuous mappings
$G^* \cc K_3$ equals the number of 1-factorizations of~$G$. 
(Consequently, there is a cut-continuous mapping $K_4 \cc K_3$, 
while there is clearly no homomorphism $K_4 \hom K_3$.)
On a similar note, let $T$ be a graph with two vertices, one edge
connecting them and one loop. It is known that the number of homomorphisms
$f: G \hom T$ equals to the number of independent sets of the graph~$G$, 
a graph parameter that is important and hard to compute. The corresponding
parameter, the number of cut-continuous mappings $g: G \cc T$ is
simple to compute (but still interesting): it is equal to the number
of cuts in~$G$, that is to $2^{|V(G)| - k}$, where $k$ is the number
of components of~$G$.

The analysis of flow problems by means of edge mappings between graphs
was pioneered by Jaeger~\cite{Jaeger};
the basic definitions were stated and developed in~\cite{DNR}.
In~\cite{NS-TT1} we studied tension-continuous (mainly
$\zet_2$-tension-continuous, that is cut-continuous) mappings more
thoroughly. Here we extend and complement results of~\cite{NS-TT1}
by treating tension-continuous mappings in an arbitrary 
abelian group instead of $\zet_2$.
We also solve several open problems from~\cite{NS-TT1}. 
Particularly, we find examples of $k$-connected graphs that are
equivalent with respect to tension-continuous mappings and not with
respect to homomorphisms (Proposition~\ref{cnhomotens} in
Section~\ref{sec:ex}).
On the positive side we give a characterization of a large class of 
graphs where tension-continuous mappings coincide with homomorphisms.
Such graphs (called here left and right homotens graphs) are studied
in Sections~\ref{sec:lh} and~\ref{sec:rh}. This also implies a shorter
proof of some results of~\cite{NS-TT1}, particularly of
universality (Theorem~\ref{niceembedding}) and
density (Theorem~\ref{AntiExtInt}) of tension-continuous mappings.
The proof of the latter uses construction~$\Delta_M(G)$
(defined in Section~\ref{sec:rh}), which is interesting in itself.

\section{Definition \& Basic Properties} \label{sec:defce}

\subsection{Basic notions---flows and tensions}

We refer to \cite{Diestel,HN} for basic notions on graphs and
their homomorphisms. 

By a graph we mean a finite directed graph with multiple
edges and loops allowed.
We write $uv$ (or sometimes $(u,v)$) for an edge from~$u$ to~$v$ (one
of them, if there are several parallel edges). A \emph{circuit} in a
graph is a connected subgraph in which each vertex is adjacent to two
edges. For a circuit~$C$, we let $C^+$ and~$C^-$ be the sets of edges
oriented in either direction. We will say that $(C^+, C^-)$ is a
\emph{splitting} of edges of~$C$.

A \emph{cycle} is an edge-disjoint union of circuits. Given a graph~$G$
and a set~$X$ of its vertices, we let~$\delta(X)$ denote the set of all
edges with one end in~$X$ and the other in~$V(G) \setminus X$; we call
each such edge set a \emph{cut} in~$G$. Let $M$ be a ring
(by this we mean an asociative ring with unity).
We say that a function $\phi: E(G) \to M$ is an \emph{$M$-flow on~$G$}
if for every vertex $v \in V(G)$  
$$
    \sum_{\hbox{$e$ enters $v$}} \phi(e) = 
    \sum_{\hbox{$e$ leaves $v$}} \phi(e) \,.
$$

A function $\tau : E(G) \to M$ is an \emph{$M$-tension
on~$G$} if for every circuit~$C$ in~$G$ (with $(C^+,C^-)$ being
the splitting of its edges) we have
$$
  \sum_{e \in C^+} \tau(e) = \sum_{e \in C^-} \tau(e) \,.
$$
We remark that for definition of flows and tensions we could use
any abelian group. But as our emphasis is on finite graphs, we are interested
in finitely generated abelian groups. Every such group is 
of form~$\zet^k \times{}\prod \zet_{n_i}^{k_i}$, therefore
we can introduce a ring structure on it. In proof of~Lemma~\ref{SIL}
we present a way how results about general abelian groups can
be inferred from finitely generated ones.

Note that $M$-tensions on a graph~$G$ form a module over~$M$
(or even a vector space, if $M$ is a field). Its dimension is $|V(G)|-k(G)$, 
where $k(G)$ denotes the number of components of~$G$. This module will
be called the \emph{$M$-tension module} of~$G$.

For a cut~$\delta(X)$ we define
$$
   \phi_X(uv) = \begin{cases}
     1  \quad \hbox{if $u \in X$ and $v \notin X$} \\
     -1  \quad \hbox{if $u \notin X$ and $v \in X$} \\
     0 \quad \hbox{otherwise}\,.\\
   \end{cases}
$$
Any such $\phi_X$ is called \emph{elementary $M$-tension}.
It is easy to prove that elementary $M$-tensions generate
the $M$-tension module.

Remark that every $M$-tension is of form $\delta p$, 
where $p:V(G) \to M$ is any mapping and 
$(\delta p) (uv) = p(v) - p(u)$ (in words, tension is a difference
of a potential).

For $M$-flows the situation is similar to $M$-tensions:
all $M$-flows on~$G$ form a module (the \emph{$M$-flow module} of~$G$)
of dimension $|E(G)| - |V(G)| + k(G)$; it is generated by \emph{elementary
flows} (those with a circuit as a support) and it is orthogonal
to the $M$-tension module.

The above are the basic notions of algebraic graph theory. For a
more thorough introduction to the subject see~\cite{Diestel}; we only mention
two more basic observations: 

A cycle can be characterized as a support of a $\zet_2$-flow and 
a cut as a support of a $\zet_2$-tension.
If $G$~is a plane graph then each cycle in~$G$ corresponds to a cut
in its dual~$G^*$; each flow on~$G$ corresponds to a tension on~$G^*$.

\subsection{Tension-continuous mappings}

The following is the principal notion of this paper: 
Let $M$ be a ring, let $G$, $G'$ be graphs and let $f: E(G) \to E(G')$ 
be a mapping between their edge sets.  We say $f$~is an
\emph{$M$-tension-continuous mapping} (shortly $TT_M$~mapping) if for
every $M$-tension~$\tau$ on~$G'$, the composed mapping $\tau f$~is an
$M$-tension on~$G$. The scheme below illustrates this definition. It
also shows that $f$ ``lifts tensions to tensions'', thus suggesting the
term $TT$~mapping. 
\begin{diagram}
      E(G) & \rTo^{f}       & E(G')     \\
           & \rdTo^{\tau f} & \dTo>\tau \\
           &                & M         \\
\end{diagram}
We write $f: G \TT[M] H$ if $f$ is a $TT_M$~mapping from~$G$
to~$H$ (or, more precisely, from~$E(G)$ to~$E(H)$).
In the important case $M = \zet_n$ we write $TT_n$ instead
of~$TT_{\zet_n}$, when $M$ is clear from the context
we omit the subscript.

Of course if $M = \zet_2$ then the orientation of edges does not
matter. Hence, if $G$, $H$ are undirected graphs and 
$f: E(G) \to E(H)$ any mapping, we say that $f$ is 
\emph{$\zet_2$-tension-continuous} ($TT_2$) if for some
(equivalently, for every) orientation $\orG$ of~$G$ and $\orH$ of~$H$, 
$f$ is $TT_2$ mapping from~$\orG$ to~$\orH$. As cuts correspond
to $\zet_2$-tensions, with this provision $TT_2$ mappings of
undirected graphs are exactly the \emph{cut-continuous} mappings:
mappings between edge sets of undirected graphs such that preimage
of every cut is a cut.

For general ring~$M$, the orientation is important. Still, we 
define that a mapping~$f: E(G) \to E(H)$ between undirected graphs $G$, $H$
is $TT_M$ if for some orientation $\orG$ of~$G$ and $\orH$ of~$H$, 
$f$ is $TT_M$ mapping from~$\orG$ to~$\orH$. This definition may
seem a bit arbitrary, but in fact it is a natural one: clearly
it is equivalent to ask that for each~$\orH$ there is an~$\orG$ such
that $f$ is a $TT_2$ mapping from~$\orG$ to~$\orH$ (we just change
orientation of edges of~$\orG$ according to change of orientation
of edges of~$\orH$). We will elaborate more on this
in Proposition~\ref{undirequiv}.

\noi {\bf Convention.} Unless specifically specified, our results
hold for both the directed and undirected case. 

Recall that $h: V(G) \to V(G')$ is called a \emph{homomorphism} if
for any $uv \in E(G)$ we have $f(u) f(v) \in E(G')$; we shortly
write $h: G \hom G'$. We define a quasiorder $\leh$ on the
class of all graphs by
$$
   G \leh G' \iff \hbox {there is a homomorphism $h: G \hom G'$.} 
$$
Homomorphisms generalize colorings: a $k$-coloring is exactly a
homomorphism $G \hom K_k$, hence $\chi(G) \le k$ iff $G \leh K_k$.
For an introduction to the theory of homomorphisms see \cite{HN}.

Motivated by the homomorphism order $\leh$, we define for 
a ring~$M$ an order $\leTT_M$ by
$$
   G \leTT_M G' \iff \hbox {there is a mapping $f: G \TT[M] G'$.} 
$$
This is indeed a quasiorder, see Lemma~\ref{compose}.
We write $G \eqTT_M H$ iff $G \leTT_M H$ and $G \geTT_M H$, 
and similarly we define~$G \eqh H$; we say $G$ and~$H$ are 
$TT_M$-equivalent, or hom-equivalent, respectively.
Occasionally, we also write $G \TT[M] H$ (instead of $G \leTT_M H$)
to denote the existence of some $TT_M$ mapping.

We define analogies of other notions used for study of homomorphisms:
a graph~$G$ is called \emph{$TT_M$-rigid} if there is no non-identical
mapping $G \TT[M] G$. Graphs $G$, $H$ are called \emph{$TT_M$-incomparable}
if there is no mapping $G \TT[M] H$, neither $H \TT[M] G$. 

If $G$ is an undirected graph, its \emph{symmetric orientation}~$\symG$ 
is a directed graph with the same set of vertices and with each
edge replaced by an oriented 2-cycle, we will 
say these two edges are opposite. The following result clarifies
the role of orientations.

\begin {proposition}   \label{undirequiv}
Let $G$, $H$ be undirected graphs, $E(H) \ne \emptyset$,
let $M$ be a ring. Then the following are equivalent.
\begin {enumerate}
  \item For some orientation $\orG$ of~$G$ and $\orH$ of~$H$ 
    it holds that $\orG \TT[M] \orH$.
  \item For each orientation $\orH$ of~$H$ exists $\orG$ of~$G$ 
    such that $\orG \TT[M] \orH$.
  \item For symmetric orientations $\symG$ of~$G$ and $\symH$ of~$H$
    it holds that $\symG \TT[M] \symH$.
\end {enumerate}
\end {proposition}

\begin {proof}
If $M = \zet_2^k$ then all statements are easily equivalent,
so suppose $M \ne \zet_2^k$.
Take a mapping $f_1: \orG \TT[M] \orH$. We may suppose
that $\orG \subseteq \symG$ and $\orH \subseteq \symH$. 
Thus if $e'$, $e''$ are opposite edges 
end $e' \in E(\orG)$, then we let $f_3(e')$ be $f_1(e')$ and 
$f_3(e'')$ be the edge opposite to $f_1(e')$. As 
cycles of $\orG$ together with the 2-cycles consisting of
opposite edges generate the cycle space of~$\symG$, 
mapping $f_3$ is $TT_M$, hence 1 implies 3.
Next take any $\orH$, suppose again $\orH \subseteq \symH$, and let
opposite edges $e'$, $e''$ of~$\symG$ correspond to $e \in G$.
At least one of the edges $f_3(e')$, $f_3(e'')$ connects the same
vertices (in the same direction) as some edge~$\bar e$ of $\orH$;
we let this one of $e'$, $e''$ to be an edge of~$\orG$ and let $f_2$ map it
to~$\bar e$.  Clearly, $f_2$ is a $TT_M$ mapping; therefore 3 implies 2. 
Finally 2 implies 1 is trivial.  
\qed
\end {proof}

\subsection{Basic properties}

In this section we summarize some properties of $TT$~mappings
which will be needed in the sequel.

\begin {lemma}   \label{compose}
Let $f: G \TT[M]H$ and $g: H \TT[M] K$ be $TT_M$~mappings. Then the
composition $g \circ f$ is a $TT_M$~mapping.
\end {lemma}

\begin {lemma}   \label{TTsubgraph}
Let $f: G \TT[M] H$, let $H'$ be a subgraph of $H$ that contains
all edges $f(e)$ for $e \in E(G)$. Then $f: G \to H'$ is $TT_M$ as well.
\end {lemma}

\begin {proof}
Take any $M$-tension $\tau'$ on~$H'$. Let $\tau' = \delta p'$
for $p': V(H') \to M$. If $V(H) = V(H')$ let $p = p'$, otherwise
extend~$p'$ arbitrarily to get~$p$. Now $\tau=\delta p$ is
an $M$-tension on~$H$ that agrees with~$\tau'$ on~$V(H')$. 
Hence $\tau' f = \tau f$, and as $\tau f$ is an $M$-tension,
$\tau' f$ is an~$M$-tension, too.
\qed
\end {proof}

An easy corollary of these observations is the
monomorphism-epimorpism factorization of $TT_M$~mappings.

\begin {corollary}   \label{factorization}
Let $f : G \TT[M] H$. Then there is a graph~$H'$ and 
$TT_M$ mappings $f_1: G \TT[M] H'$, $f_2: H' \TT[M] G$
such that $f_1$~is surjective and $f_2$~injective.
\end {corollary}

Another easy (but useful) way to modify $TT_M$~mapping is 
by adding parallel edges. The next result shows, that we may
in many respects restrict ourselves to bijective $TT_M$~mappings
(this approach was taken by~\cite{LMT,Kelmans}).
A bijection $G \TT[M] H$ may be viewed as an identification
$E(G) = E(H)$, therefore we in fact study when
the tension module of~$H$ is a submodule of tension module of~$G$
(this language was used in~\cite{Shih-thesis}).

\begin {lemma}   \label{parallel}
Let $f: G \TT[M] H$ be a $TT_M$~mapping of (directed or undirected)
graphs. Then there is a graph~$H'$ and a mapping $f': E(G) \to E(H')$ 
such that
\begin {itemize}
  \item $f'$~is $TT_M$,
  \item $f'$ is bijective,
  \item we can get $H'$ by adding parallel edges and 
     deleting edges from~$H$.
  \item for each edge $a \in E(G)$ the edge $f'(a)$ connects
     the same vertices as $f(a)$. 
\end {itemize}
\end {lemma}

\begin {proof}
For an edge $e \in E(H)$ we let $c(e) = |f^{-1}(e)|$ be the 
number of edges that map to~$e$. We replace each edge of~$H$
by~$c(e)$ parallel edges in the same direction (in case of directed
graphs) as~$e$ and keep all vertices; we let $H'$ denote the resulting graph.
We define $f'(a)$ to be any one of the parallel edges that
replaced~$f(a)$, making sure that $f'$ is injective (therefore 
bijective). Clearly, for any $p: V(H) = V(H') \to M$, if 
we consider the $M$-tensions $\tau = \delta p$ of $H$ 
and $\tau' = \delta p$ of $H'$, then $f\circ \tau = f'\circ \tau'$.
Thus if $f$ was a $TT_M$ mapping, $f'$ is $TT_M$ as well.
\qed
\end {proof}

If $C$ is a circuit with a splitting~$(C^+, C^-)$, we say that 
$C$~is \emph{$M$-balanced} if $(|C^+|-|C^-|)\cdot 1 = 0$
(with $0$, $1$, and operations in~$M$).
Otherwise, we say $C$~is \emph{$M$-unbalanced}.
Let $g_{\sss M}(G)$ denote the length of the shortest $M$-unbalanced
circuit in~$G$, if there is none we put $g_M(G) = \infty$.
For the particular case $M = \zet_2$, a circuit is $M$-balanced
if it is even, hence $g_{\sss \zet_2}(G)$ is the odd-girth of~$G$.
We also have $G \TT[M] \edge$ iff any constant mapping $E(G) \to M$ is an
$M$-tension. This clearly happens precisely when all circuits in~$G$
are $M$-balanced, equivalently, if $g_{\sss M}(G) = \infty$.
As a consequence of this, the function $g_{\sss M}$ provides us with an
invariant for the existence of $TT_M$~mappings, as shown in the next
two lemmas.

\begin {lemma}   \label{BalancedCircuit}
Let $M$ be a ring, let $G$, $H$ be directed graphs, let
$f : G \TT[M] H$. 
If $C$ is an $M$-unbalanced circuit in~$G$
then $f(C)$ contains an $M$-unbalanced circuit.
\end {lemma}

\begin {proof}
The inclusion homomorphism $C \to G$ induces a $TT_M$ mapping,
composition with~$f$ yields $C \TT[M] H$. By Lemma~\ref{TTsubgraph} we
get a mapping $C \TT[M] f(C)$. If all circuits in~$f(C)$ are
$M$-balanced, then $f(C) \TT[M] \edge$ and, by composition we have 
$C \TT[M] \edge$.  This contradicts the fact that $C$ is $M$-unbalanced.
\qed
\end {proof}

\begin {lemma}   \label{invariant}
Let $G \leTT_M H$ be directed graphs.
Then $g_{\sss M}(G) \ge g_{\sss M}(H)$. 
\end {lemma}

\begin {proof}
If $g_{\sss M}(G)=\infty$, the conclusion holds. Otherwise, let $C$
be an $M$-unbalanced circuit of length~$g_{\sss M}(G)$ in~$G$. 
By Lemma~\ref{BalancedCircuit}, $f(C)$ contains an $M$-unbalanced
circuit. It is of size at least~$g_{\sss M}(H)$ and at most~$g_{\sss M}(G)$.
\qed
\end {proof}

An alternative definition of tension-continuous mappings 
(proved in~\cite{DNR}) is often useful.
For mappings $f: E(G) \to E(H)$ and $\phi: E(G) \to M$ we let $\phi_f$
denote the \emph{algebraical image of $\phi$}: that is we define a mapping
$\phi_f:E(H) \to M$ by 
$$
  \phi_f(e') = \sum_{e \in f^{-1}(e')} \phi(e) \,.
$$


\begin {lemma}   \label{altdef}
Let $f: E(G) \to E(H)$ be a mapping. Then $f$ is
$M$-tension-continuous if and only if for every
$M$-flow $\phi$ on~$G$, its algebraical image $\phi_f$
is an~$M$-flow. Moreover, it is enough to verify this
property for the basis of the flow module (elementary
flows supported by an elementary cycle).

We formulate this explicitly for $M=\zet_2$. Mapping $f$ is cut-continuous 
if and only if for every cycle~$C$ in~$G$, the set of edges
of~$H$, to which an odd number of edges of~$C$ maps, is a cycle.
\end {lemma}

For a homomorphism (of directed or undirected graphs) 
$h: V(G) \to V(G')$ we let~$h^\sharp$ denote the 
\emph{induced mapping on edges}, that is 
$h^\sharp ((u,v)) = (h(u),h(v))$, or 
$h^\sharp (\{u,v\}) = \{h(u),h(v)\}$.
If $h$ is an \emph{antihomomorphism}, that is for every edge 
$(u,v)\in E(G)$ we have $(h(v),h(u)) \in E(G')$ ($h$ reverses every edge), 
we define $h^\flat((u,v)) = (h(v),h(u))$ and call it a mapping
induced by antihomomorphism.
If $G'$ has parallel edges, then $h^\sharp$ is not unique: we just
ask that $h^\sharp$ maps each of the edges $(u,v)$ to some of the
edges $(h(u),h(v))$; similarly for homomorphisms of undirected graphs
and for antihomomorphisms.
The following easy lemma is the starting point of our investigation.

\begin {lemma}   \label{homo}
Let $G$, $H$ be (directed or undirected) graphs, $M$ a ring.
For every (anti)ho\-momorphism $f$ from~$G$ to~$H$ the induced mapping
$f^\sharp$ ($f^\flat$, respectively) from~$G$ to~$H$
is $M$-tension-continuous.
Consequently, from $G \leh H$ follows $G \leTT_M H$.
\end {lemma}

\begin {proof}
It is enough to prove Lemma~\ref{homo} for homomorphisms of directed 
graphs. So let $f: G \to H$ be such homomorphism, 
$\phi: V(H) \to M$ a tension. We may
assume that $\phi$ is an elementary tension corresponding to the cut 
$\delta(X)$. Then the cut $\delta(f^{-1}(X))$
determines precisely the tension $\phi\circ f$.
\qed
\end {proof}

The main theme of this paper is to find similarities and differences
between orders $\leh$ and $\leTT_M$. 
In particular we are interested in when the converse to Lemma~\ref{homo}
holds. Now, we present a more precise version of
Question~\ref{q:basic} stated in the introduction.

\begin{problem} \label{homoTT}
Let $f : E(G) \to E(H)$. Find suitable conditions for $f$, $G$, $H$ 
that will guarantee that whenever $f$ is $TT_M$, then it
is induced by a homomorphism (or an antihomomorphism); 
i.e.\ that there is a homomorphism (or an antihomomorphism) 
$g : V(G) \to V(H)$ such that $f = g^\sharp$ (or $f = g^\flat$).
\end{problem}

Shortly, we say a mapping is \emph{induced}
if it is induced by a homomorphism or an antihomomorphism.
Problem~\ref{homoTT} leads us to the following definitions.

\begin{definition} \label{def:homotens}
We say a graph~$G$ is \emph{left $M$-homotens} if 
for every loopless graph~$H$ every $TT_M$~mapping from~$G$ to~$H$
is \emph{induced} (that is induced by a homomorphism 
or an antihomomorphism). For brevity we will often
call left $M$-homotens graphs just $M$-homotens graphs
(following~\cite{NS-TT1}).

On the other hand, $H$ is a \emph{right $M$-homotens} graph if for every
graph~$G$ statements $G \hom H$ and $G \TT[M] H$ are equivalent. 
\end{definition}

We should note here, that the precise analogy of left $M$-homotens 
graphs---every $TT_M$ mapping is induced---is not interesting, as
this is much too strong requirement. For simplicity, suppose $M=\zet_2$.
Let $H$ be such graph, let $\Delta(H)$ be as defined before
Lemma~\ref{cube}. The mapping $f:\Delta(H) \TT[2] H$ given by
$f(\{A,B\}) = A \sym B$ is induced by an (anti)homomorphism, say~$g$.  
Now this can happen only if for every $A \in V(\Delta(H))$ vertex
$g(A)$ is adjacent to every edge $e$ of $H$.  (To see this, note that
$f(\{A, A \sym e\}) = e$, therefore $g(A)$ is one of the end vertices of~$e$.)
And this in turn can happen only if $H$ is edgeless, or if $H = K_2$. 

Definition of left $M$-homotens makes sense for both directed
and undirected graphs. 
If $M=\zet_2^k$ then there are only trivial directed $M$-homotens graphs
(namely an orientation of a matching). Thus, we restrict to study
of undirected homotens graphs in this case~\cite{NS-TT1}.
For other rings, Proposition~\ref{reorhomotens} states that
the orientation does not play any role; this will be 
useful in Section~\ref{sec:lh} in our study of directed $M$-homotens graphs.

For $M \ne \zet_2^k$ we might study undirected $M$-homotens graphs,
too. The relationship between these two notions (undirected
graph is homotens versus some its orientation is homotens) is 
not clear. For every~$M$, the latter notion implies the former one;
however, somewhat surprisingly, both notions
are equivalent for many rings~$M$ (at least for such, in which
the equation~$x+x=0$ has no nonzero solution). 
(For right homotens graphs, the above discussion applies, too.)

\begin {proposition}   \label{reorhomotens}
Let $G_1$, $G_2$ be two directed graphs, such that we can get
$G_2$ from~$G_1$ by changing directions of edges, deleting
and adding multiple edges. Let $M$ be a ring.
Then $G_1$ is left $M$-homotens if and only if
$G_2$ is left $M$-homotens.
\end {proposition}

\begin {proof}
Suppose $G_1$ is not homotens, that is there is a graph~$H_1$
and a mapping~$f_1:G_1 \TT[M] H_1$ that is not induced. 
By Lemma~\ref{parallel} we may suppose that $f_1$~is injective.
We modify $f_1$ and $H_1$, to get a non-induced
mapping~$f_2: G_2 \TT[M] H_2$.
If we change an orientation of an edge, we change an orientation
of the corresponding edge in~$H_1$. If we add an edge parallel
to some edge $e$ of~$G_1$ then we map it to a new edge of~$H_1$, 
parallel to~$f_1(e)$. It is clear, that we get a~$TT_M$~mapping
that is not induced.
\qed
\end {proof}

\section{Examples} \label{sec:ex}

We illustrate the complex relationship of homomorphisms and $TT$~mappings
by several examples presenting the similarities and (mainly) the differences
in concrete independent settings. 
Towards the former, we provide an infinite chain
and antichain of~$\leTT_{\sss \zet_2}$, 
thereby exhibiting a similar behaviour
of homomorphisms and $TT$~mappings. On the other hand, we 
show that arbitrarily high connectivity of the source and target graphs
does not force $TT_\zet$~mappings (much the less $TT_M$~mappings) and
homomorphisms to coincide.  Finally, we show that an equivalence class
of~$\eqTT_{\sss \zet_2}$ can contain exponentially many equivalence
classes of~$\eqh$.

Proposition~\ref{illustration} appears already in~\cite{DNR},
we include a proof for the convenience of the reader.
Note that this result will be strongly generalized 
by Theorems~\ref{homotens},~\ref{niceembedding}, and~\ref{AntiExtInt}. 

\begin {proposition}   \label{illustration}
Graphs $K_{2^t}$ form a strictly increasing chain 
in $\leTT_{\sss \zet_2}$ order, that is 
$K_4 \lTT_{\sss \zet_2} K_8 \lTT_{\sss \zet_2} K_{16} \lTT_{\sss \zet_2} \cdots$.
There are graphs $G_1$, $G_2$, \dots\ that form an infinite antichain:
there is no mapping $G_i \TT[2] G_j$ for $i \ne j$.
\end {proposition}

\begin {proof}
By Proposition~6 of~\cite{DNR} (compare also Corollary~\ref{rhcomplete}
of this paper), for any graph $G$
\begin{equation}\label{eq:equiv}
  G \hom K_{2^k}  \Longleftrightarrow G \TT[2] K_{2^k} \,.
\end{equation}
This implies the first
part. For the second part, let $G_t$ be the Kneser graph
$K(n,k)$ with $k=t(2^t-2)$ and $n=2k+2^t-2$. It is known 
that $\chi(G_t) = n -2k + 2 = 2^t$. This by equivalence $(\ref{eq:equiv})$
implies that $G_i \nTT[2] G_j$ for $i > j$. 
The remaining part follows from Lemma~\ref{invariant}: 
It is known that the shortest odd cycle in $K(n,k)$ is
the smallest odd number greater or equal to~$n/(n-2k)$, which
means that $g_{\sss \zet_2}(G_t) = 2t+1$. 
\qed
\end {proof}

The differences of $TT$~mappings and homomorphisms are easy to
find. For example let $\{e_1,e_2,e_3\}$ be the edges of~$K_3$, and
color the edges of~$K_4$ properly by three colors. 
We send both edges of color~$i$ to~$e_i$. This mapping is easily checked
to be $TT_2$, so we have $K_4 \TT[2] K_3$ but obviously there is no
homomorphism $K_4 \to K_3$.
On the contrary, $TT_\zet$~mappings are more restricted and, indeed, there is
no $TT_\zet$ mapping from an orientation of~$K_4$ to an orientation of~$K_3$.
A simple example of $TT_\zet$ mapping that is not induced by a
homomorphism is a noncyclic permutation of edges of an oriented
circuit. E.g., let $E(\orC_5) = \{e_0, e_1, \dots, e_4\}$ in this
order, and define $f(e_i) = e_{2i \bmod 5}$. Then $f$ is $TT_\zet$, 
on the other hand, $f$~maps adjacent edges to nonadjacent edges, hence
is not induced by a homomorphism. By applying the arrow
construction---that is by replacing each oriented edge by a suitable
graph (see \cite{HN} and also proof of Proposition~\ref{manydif} for more
details) it is easy to produce graphs $G$, $H$ such that $G \TT[\zet] H$
but $G \nhom H$. No graphs $G$, $H$ obtained in this manner are 3-connected;
Whitney's theorem (two 3-regular graphs with the same cycle matroid are
isomorphic) seems to suggest, that this situation may not repeat for graphs
with higher connectivity. 
Therefore, the following lemma may be a bit surprising.

\begin {proposition}   \label{cnhomotens}
For every $k$ there are $k$-connected graphs $G$, $H$ such that $G
\TT[\zet] H$ but $G \nhom H$.  Therefore, for each~$k$ exists a
$k$-connected graph that is not $\zet$-homotens. 
\end {proposition}

\begin {figure}
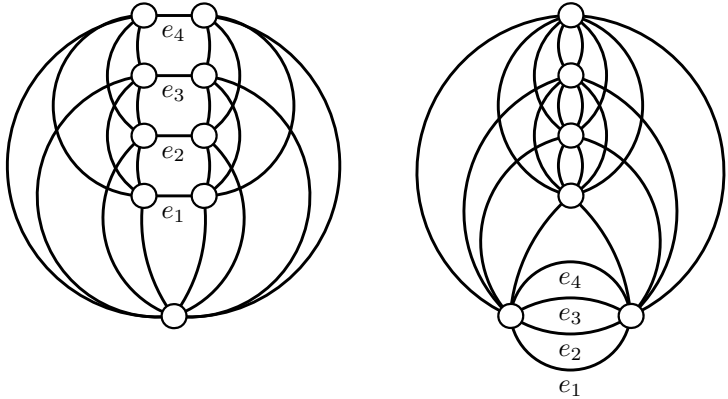

\centerline {
  \vtop{\kern 0pt\hbox{\includegraphics{cnhomotens.0}}}
  \hfil
  \vtop{\kern 0pt\hbox{\includegraphics{cnhomotens.1}}}
}
  \caption{The left graph is an example of highly connected graph that is not
    $\zet$-homotens; the right one is a witness for the former 
    not being $\zet$-homotens.}
  \label{fig:cnhomotens}
\end {figure}

\begin {proof}
Fix a $k$, let $G$, $H$ be graphs illustrated for $k=4$ 
in Figure~\ref{fig:cnhomotens}.
\footnote{If we wish to construct directed graphs, consider
any orientation of them, such that corresponding edges of~$G$ and
of~$H$ are oriented in the same way.}
(The construction is due to Shih~\cite{Shih-thesis}.)

Clearly both $G$ and $H$ are $k$-connected and there 
there is no homomorphism between them.
The natural bijection between $G$ and $H$---we identify 
the left $K_k$'s in $G$ and~$H$, the right $K_k$'s in $G$ and~$H$, 
and the edges $e_i$ as depicted in the Figure---is easily
checked to be $TT_\zet$. 
\qed
\end {proof}

Further examples of graphs with negative answer to Problem~\ref{homoTT}
are listed in~\cite{NS-TT1}, here we only mention the perhaps
most spectacular example: Petersen graph admits a $TT_2$~mapping
to~$C_5$. This mapping (and many others) may be obtained using
the following construction: 
Given an (undirected) graph $G = (V,E)$ write $\Delta(G)$
for the graph $({\cal P}(V),E')$, where $AB \in E'$ iff
$A \Delta B \in E$ (here ${\cal P}(V)$ denotes the set of all
subsets of $V$ and $A \Delta B$ the symmetric difference of sets
$A$ and $B$).

\begin{lemma} \label{cube} 
Let $G$, $H$ be undirected graphs. 
Then $G \TT[2] H$ iff $G \hom \Delta(H)$.
\end{lemma}

We can formulate analogous construction and result for rings
$M \ne \zet_2$; this is done in Section~\ref{sec:delta}.
We conclude this section by a more quantitative example. 

\begin {proposition}   \label{manydif}
There are $2^{cn}$ undirected graphs with $n$ vertices that form an
antichain in the homomorphism order, yet all of them are
$TT_2$-equivalent. 
\end {proposition}

\begin {figure}
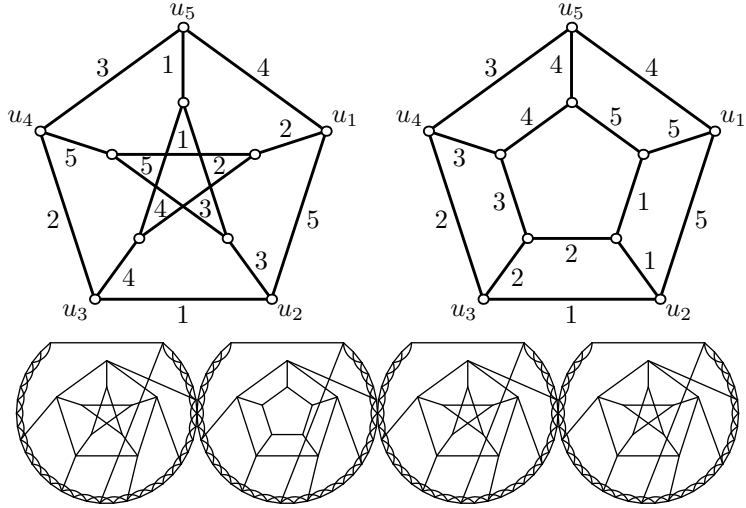

\centerline{
  \hfil
  \includegraphics{Pet.0}
  \hfil
  \includegraphics{Pet.1}
  \hfil
}
\medskip
\centerline{
  \includegraphics{niceindicator.1}
}
  \label{fig:manydif}
  \caption{Petersen graph and the prism of~$C_5$---two $TT_2$-equivalent
    graphs used in the proof of Proposition~\ref{manydif}. Below is an example 
    of the construction for $n=4$, $t=(1,0,1,1)$.}
\end {figure}

\begin {proof}
To simplify notation, we will construct $\binom n{\floor{n/2}}$ 
graphs with $sn+1$~vertices, this clearly proves the proposition.
We use the {\it replacement operation} of~\cite{HN}. Let $H$ be a
graph (we explain later how do we choose it), let $a$, $b$, $x_1$,
\dots, $x_5$ be pairwise distinct vertices of~$H$. Next, we take an
oriented path with $n$ edges and replace each of them by a copy
of~$H$. That is, we take $H_1$, \dots, $H_n$---isomorphic copies
of~$H$---and identify vertex $b$ of~$H_i$ with $a$ of $H_{i+1}$ 
(for every $i = 1, \dots, n-1$). Let $G$~be
the resulting graph. 

Finally, for each $t \in \{0,1\}^n$ we present a graph~$G_t$.
We let $F_i$ be a copy of the Petersen graph~$P$ if $t_i = 1$, 
and a copy of the prism of~$C_5$---graph $R$ in
Figure~\ref{fig:manydif}---if $t_i=0$.
We construct the graph~$G_t$ as a vertex-disjoint union 
of~$G$, $F_1$, \dots, $F_t$ plus some `connecting edges': 
for every $i=1, \dots, n$ and $j=1, \dots, 5$ we let $x_i^j$ 
denote the copy of~$x_j$ in~$H_i \subset G$ and
$u_i^j$ the copy of $u_j$ in~$F_i$; we let $x_i^j u_i^j$ 
be an edge of~$G_t$. 
Note that each $G_t$ has $(|V(P)| + |V(H)| - 1)n + 1$ vertices.

\textbf{Claim 1.} 
\begin{minipage}[t]{0.8\textwidth}
$H$ can be chosen so that the only homomorphism 
$G \to G$ is the identity. Moreover the vertices $x_i$ can be
chosen so that the distance between any two of them is at least~$4$.
\end{minipage}

This follows immediately from techniques of~\cite{HN}, 
e.g. we can take $H_9$ from the Figure~4.9 of~\cite{HN} as our graph~$H$.

\textbf{Claim 2.} If $G_t \hom G_{t'}$ then $t_i \le t'_i$ holds for each $i$.

Take any homomorphism $f: G_t \hom G_{t'}$, fix an~$i$, and
let $F_i$ ($F'_i$) be the copy of~$P$ or~$R$ that 
constitute the $i$-th part of graph~$G_t$ ($G_{t'}$ respectively). 
By Claim~1, $f$~maps the vertices of~$G$ identically, in particular
$f(x_i^j) = x_i^j$. As the only path of length~$3$ connecting
vertices $x_i^j$ and $x_i^{j \bmod 5 + 1}$ is the one containing 
vertices $u_i^j$ and $u_i^{j \bmod 5 + 1}$, mapping $f$ satisfies 
$f(u_i^j) = u_i^j$ as well. Consequently, $f$ maps vertices
of~$F_i$ to vertices of~$F'_i$. To show $t_i \le t'_i$
it remains to observe that there is no homomorphism $P \hom R$.

\textbf{Claim 3.} For every $t$, $t'$ we have $G_t \TT[2] G_{t'}$. 

We map every edge of $G$ and every edge $x_i^j u_i^j$ and 
$u_i^j u_i^{j \bmod 5 + 1}$ identically (we call such edges
\emph{easy edges}). 
We map edges of $F_i$ in~$G_t$ to edges of the outer pentagon of $F_i$
in~$G_{t'}$ by sending an edge to the outer edge with the same number in
Figure~\ref{fig:manydif}. 
To check that this is indeed a $TT_2$~mapping we use
Lemma~\ref{altdef}: if $C$ is a cycle contained in some~$F_i$ then we
easily check that algebraical image of~$C$ is a cycle. If $C$ 
contains only easy edges that it is mapped identically, so its
algebraical image is again a cycle. As every cycle can be written
as a symmetric difference of these two types, we conclude that 
we have constructed a $TT_2$~mapping.

Now we are ready to finish the proof. Consider a set~$A$ containing
all vertices of~$\{0,1\}^n$ with $\floor{n/2}$ coordinates equal to~1.
By Claim~2, graphs $G_t$, $G_{t'}$ are homomorphically incomparable 
for distinct $t, t' \in A$.
On the other hand, by Claim~3, all of the graphs are
$TT_2$-equivalent.
\qed
\end {proof}

In this proof we can use other building blocks instead of 
Petersen graph and the pentagonal prism. To be concrete, we can
take graphs $G$, $H$ from Proposition~\ref{cnhomotens} and 
use graphs $G \dunion H$ and $H \dunion H$. If we slightly modify the
construction, we can prove version of Proposition~\ref{manydif} for
$TT_\zet$ mappings, and therefore for $TT_M$ mappings for
arbitrary~$M$.  Moreover, by another small change of the construction,
we can guarantee that all of the constructed graphs are $k$-connected
(for any given~$k$).

It would be interesting to know if $2^{cn}$ from Proposition~\ref{manydif}
can be improved. Note that in the homomorphism order $\leh$
the maximal antichain has full cardinality \cite{KR},
that is there are
$$
  \frac{1}{n!} \binom { \binom n2 }{ \floor{\frac 12 \binom n2}} (1-o(1))
$$
homomorphically incomparable graphs with $n$-vertices.
Proposition~\ref{manydif} claims that at least $2^{cn}$ of these
graphs are contained in one equivalence class of~$\eqTT_M$.
  
\section{Left homotens graphs} \label{sec:lh}

In this section we point out similarities between homomorphisms
and $TT_M$~mappings by defining a class of graphs that
force any $TT_M$ mapping from them to be induced. 
We prove a surprising result that most graphs have this property. 
In Section~\ref{sec:applications} we use these graphs to 
find an embedding of category of graphs and homomorphism to 
the category of graphs and $TT_M$ mappings,
simplifying and generalizing a result of~\cite{NS-TT1}.
 
\subsection{A sufficient condition}
\label{sec:sufficient}

Recall (Definition~\ref{def:homotens}) that a graph~$G$ is left
$M$-homotens if every $TT_M$~mapping from~$G$ (to any graph) is
induced. The characterization 
of left $M$-homotens graphs seems to be a difficult problem; 
in this section we obtain a general sufficient condition in terms
of~\emph{nice} graphs. This notion was introduced
and proved to be a sufficient condition in~\cite{NS-TT1}
but only for~$M=\zet_2^k$. 
Here, we prove it to be sufficient for all rings
\emph{different} from~$\zet_2^k$. 
(Restricting to $M \ne \zet_2^k$ enables us to slightly weaken the
sufficient condition.)

In Proposition~\ref{cnhomotens} we saw that high connectivity 
does not imply homotens. In Corollary~\ref{triangle} we will see
that every vertex of a homotens graph is incident with a triangle.
In view of this, a sufficient condition for homotens has to be
somewhat restrictive.

\begin {definition}   \label{def:wnice}
We say that an undirected graph~$G$ is \emph{nice} if the following holds
\begin {enumerate}
  \item every edge of~$G$ is contained in some triangle
  \item every triangle in~$G$ is contained in some copy of~$K_4$
  \item every copy of~$K_4$ in~$G$ is contained in some copy of~$K_5$
  \item for every $K$, $K'$ that are copies of~$K_4$ in~$G$ there is a 
    sequence of vertices
    $v_1$, $v_2$, \dots, $v_t$ such that
    \begin {itemize}
      \item $V(K) = \{v_1, v_2, v_3, v_4\}$,
      \item $V(K') = \{v_{t}, v_{t-1}, v_{t-2}, v_{t-3}\}$, 
      \item $v_i v_j$ or $v_j v_i$ is an edge of~$G$ whenever
        $1 \le i < j \le t$ and $j \le i+3$.
    \end {itemize}
\end {enumerate}

We say that a graph is \emph{weakly nice} if conditions 1, 2, and 4
in the list above are satisfied.
Finally, we say that a \emph{directed} graph is (weakly) nice, 
if the underlying undirected graph is (weakly) nice.
\end {definition}

Before we prove Theorem~\ref{wniceTT}, which we are aiming to, 
we restate here analogous result that appears as Theorem~13
in~\cite{NS-TT1}.

\begin {theorem}   \label{niceTT}
Let $G$, $H$ be undirected graphs, let $G$ be nice, 
and let $f: G \TT[2] H$. 
Then $f$ is induced by a homomorphism of the underlying undirected
graphs. Shortly, every undirected nice graph is $\zet_2$-homotens.  
\end {theorem}

\begin {theorem}   \label{wniceTT}
Let $G$, $H$ be (directed or undirected) graphs, let $G$ be weakly nice, 
let $M\ne (\zet_2)^r$ any ring. Suppose $f : G \TT[M] H$. 
Then $f$ is induced by a homomorphism or an 
antihomomorphism. Shortly, every weakly nice graph is $M$-homotens.
\end {theorem}

We take time out for a technical lemma.

\begin {lemma}   \label{TTK4}
Let $M$ be a ring that is not isomorphic to a power 
of~$\zet_2$.
Let $f: \orK_4 \TT[M] H$, where $H$ is any loopless graph and
$\orK_4$ any orientation of~$K_4$.
Then $f$ is induced by an injective homomorphism or antihomomorphism.
Moreover, this (anti)ho\-momorphism is uniquely determined.
\end {lemma}

\begin {proof}
Suppose first that $f(\orK_4)$ is a three-colorable graph, i.e., that  
there is a homomorphism $h: f(\orK_4) \to \symK_3$, where $\symK_3$ is
the directed graph with three vertices and all six oriented
edges among them.  
A composition of $TT_M$~mapping $f : \orK_4 \TT[M] f(\orK_4)$
with~$h^\sharp$ gives $g:\orK_4 \TT[M] \symK_3$. 
Consider the three cuts of size $4$ in $\symK_3$: $X_1$, $X_2$, $X_3$.
As $M$ is not a power of~$\zet_2$, $1 + 1 \ne 0$; 
let $\phi_i$ be $M$-tension that attains value $\pm 1$ on~$X_i$
and $0$ elsewhere. We can choose $\phi_i$ so, that for every 
$e \in E(\symK_3)$ we have $\{\phi_1(e), \phi_2(e), \phi_3(e)\} = \{0, \pm 1\}$. 
As $g$ is $TT_M$, mappings $\psi_i = \phi_i g$ are $M$-tensions and for every 
$e \in E(\orK_4)$ we have $\{ \psi_1(e),\psi_2(e), \psi_3(e)\} = \{0, \pm 1\}$.
(*)

Call an $M$-tension \emph{simple} if it attains only values~$0$ 
and~$\pm 1$. We will show that three simple
$M$-tensions $\psi_1$, $\psi_2$, $\psi_3$ on~$\orK_4$ with property (*)
do not exist.

To this end, we will characterize sets $\Ker \psi = \{e\in E(\orK_4),
\psi(e) = 0\}$ for simple $M$-tensions $\psi$. Let $\psi$
be such tension. Pick $v \in V(\orK_4)$ and let $e_1$, $e_2$, $e_3$ be
adjacent to~$v$. 
Note that $\psi$ is determined by its values on $e_1$, $e_2$, $e_3$.
We may suppose that each $e_i$ is going out of~$v$;
otherwise we change orientation of some edges and the sign of $\psi$ on them.
Further, we may suppose that 
$|\{i, \psi(e_i) = 1\}| \ge |\{i, \psi(e_i) = -1\}|$;
otherwise we consider $-\psi$.
Thus, we distinguish the following cases (see Figure~\ref{fig:proofTTK4}).
\begin {itemize}
  \item $\psi(e_i) \in \{0, 1\}$ for each $i$. \\
    Let $z$ be the number of $e_i$ such that $\psi(e_i) = 0$. Then
    $\psi$ is generated by a cut with $z+1$ vertices on one side of
    the cut. Therefore, the set $\Ker \psi$ is either the edge set of
    a~$\orK_4$, of a triangle, or it is a pair of disjoint edges.
  \item $\psi(e_1) = 1$, $\psi(e_2) = 0$, $\psi(e_3) = -1$. \\
    In this case $\Ker \psi$ is a single edge.
    Note, that this case (and the next one) may happen only if 
    $1 + 1 + 1 = 0$.  
  \item $\psi(e_1) = \psi(e_2) = 1$, $\psi(e_3) = -1$. \\
    In this case too, $\Ker \psi$ is a single edge.
\end {itemize}

\begin{figure}
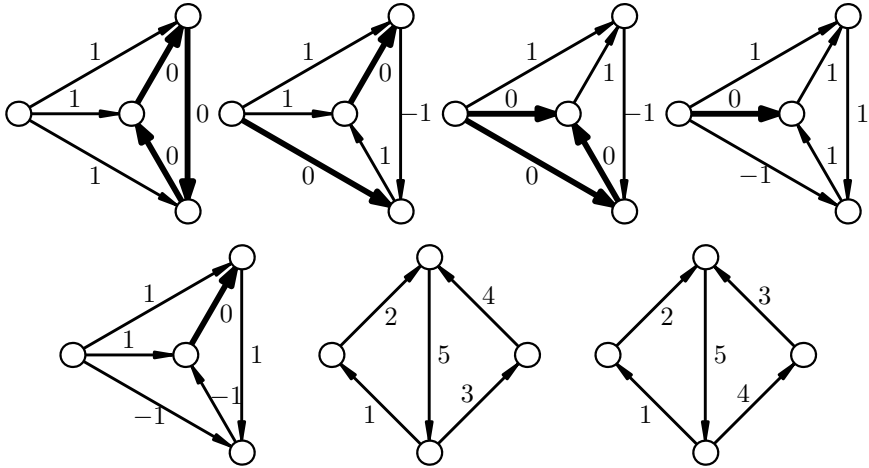

\centerline{
  \includegraphics{proofnice.8}
  \hfil
  \includegraphics{proofnice.9}
  \hfil
  \includegraphics{proofnice.10}
  \hfil
  \includegraphics{proofnice.11}
}
\medskip
\centerline{
  \includegraphics{proofnice.12}
  \hfil
  \includegraphics{proofnice.13}
  \hfil
  \includegraphics{proofnice.14}
}
  \caption{Illustration of proof of Lemma~\ref{TTK4}.}
  \label{fig:proofTTK4}
\end{figure}

Hence, $E(\orK_4)$ is partitioned into three sets, whose sizes
are in~$\{1, 2, 3, 6\}$. Therefore, there are two possibilities:
\begin {itemize}
  \item $6 = 3 + 2 + 1$: The complement of a triangle is a star of three
    edges, there are no two disjoint edges in it.
  \item $6 = 2 + 2 + 2$: In this case, all three $\psi_i$'s are
    generated by a cut. Suppose $\orK_4$ is oriented as in
    Figure~\ref{fig:proofTTK4}, the values of~$\psi_1$ are indicated.
    It is not possible to fulfill the condition (*) on both 
    edges from $\Ker \psi_1$.
\end {itemize}

So far we have proved, that the chromatic number of $f(\orK_4)$ is at
least four. As $f(\orK_4)$ has at most 6 edges, its chromatic number is exactly
four. Let $V_1$, \dots, $V_4$ be the color classes. There is exactly
one edge between two distinct color classes (otherwise
the graph is three-colorable). Thus, $f$ is a bijection.
Next, $|V_i| = 1$ for every~$i$ 
(as otherwise, we can split one color-class to several 
pieces and join these to the other classes; again, the
graph would be three-colorable). Consequently, $f(\orK_4)$
is some orientation of~$K_4$.

We call \emph{star} a set of edges sharing a vertex. If we let $\phi$ be
a simple $M$-tension on~$f(\orK_4)$ corresponding to a cut which is a
star, then $\phi f$ is a simple tension that is nonzero exactly on
three edges ($f$ is a bijection). By the characterization of zero sets of 
simple tensions we see that preimage of each star is a star.
As $f$ is a bijection and preimage of every star is a star, also image of
every star is a star. This allows us to define a vertex bijection
$g: V(\orK_4) \to V(f(\orK_4))$ by letting $g(u) = u'$ iff
the $f$-image of the star with~$u$ as the central vertex
is the star centered at~$u'$.
Stars sharing an edge map to stars sharing an
edge, hence $f$ is induced by~$g$, which is either a homomorphism 
or an antihomomorphism.
\qed
\end {proof}

\begin {proof}[Theorem~\ref{wniceTT}]
It is convenient to suppose that $G$~contains no 
parallel edges (Proposition~\ref{reorhomotens}).  Let $K$ be a copy of~$K_4$
in~$G$ (by this we mean here that $K$~is some orientation of~$K_4$).
By Lemma~\ref{TTK4} the restriction of~$f$ to~$K$ is induced by an
(anti)homomorphism, let it be denoted by~$h_K$. That is, we assume
$f|_{E(K)} = h_K^\sharp$ (or $f|_{E(K)} = h_K^\flat$).

As every edge is contained in some copy of~$K_4$, 
it is enough to prove that there is a common extension
of all mappings $\{h_K \mid K \subseteq G,\ K \simeq K_4\}$
(we may define it arbitrarily on isolated vertices of~$G$).

We say that $h_K$ and~$h_{K'}$ \emph{agree} if for
any $v \in V(K) \cap V(K')$ we have $h_K(v) = h_{K'}(v)$ and 
either both $h_K$, $h_{K'}$ are homomorphisms or 
both are antihomomorphisms.
Thus, we need to show that any two mappings $h_K$, $h_{K'}$ agree.

First, let $K$, $K'$ be copies of $K_4$ that intersect in a triangle. 
Then $h_K$ and~$h_{K'}$ agree (note that this does not necessarily
hold if the intersection is just an edge, see
Figure~\ref{fig:proofTTK4}). 

Now suppose $K$, $K'$ are copies of $K_4$ that have a common vertex~$v$. 
Since the graph~$G$ is weakly nice, we find $v_1$, $v_2$, \dots, $v_t$ as 
in Definition~\ref{def:wnice}.
Let $G_i = G[\{v_i, v_{i+1}, v_{i+2}, v_{i+3}\}]$: every 
$G_i$ is a copy of $K_4$, $G_1 = K$ and $G_{t-3}=K'$.
Suppose $v = v_l = v_r$, where $l \in \{1, 2, 3, 4\}$, 
$r \in \{t-3, t-2, t-1, t\}$. 
Consider a closed walk $W = v_{l}, v_{l+1}, \dots, v_{r-1}, v_r$.
Let $v'_i = h_{G_i}(v_i)$ for $l \le i \le r-3$
and $v'_i = h_{G_{r-3}}(v_i)$ for $r-3 \le i \le r$. 
Mappings $h_{G_i}$ and $h_{G_{i+1}}$ agree, hence
$v'_i v'_{i+1} = f(v_i v_{i+1})$ is an edge of~$H$. So
$W' = v'_{l}, v'_{l+1}, \dots, v'_{r-1}, v'_r$ is a walk in $H$.

Let $\phi$ be `a $\pm 1$-flow around~$W$', formally
$$
  \phi(e) = \sum_
            {\begin{array}{c}
               \scriptstyle l \le i \le r-1 \\
               \scriptstyle e = (v_i,v_{i+1})
              \end{array}
            } 1
         - \sum_
            {\begin{array}{c}
               \scriptstyle l \le i \le r-1 \\
               \scriptstyle e = (v_{i+1},v_i)
              \end{array}
            } 1 \,.
$$
Clearly $\phi$ is an $M$-flow. Similarly, define $\phi'(e)$ from~$W'$. 
We have $\phi' = \phi_f$, hence $\phi'$ is a flow (Lemma~\ref{altdef}).
This can happen only if $W'$ is a closed walk, 
that is $v'_l = v'_r$. 

By definition, $v'_r = h_{K'}(v)$. As mappings $h_{G_i}$ and $h_{G_{i+1}}$
agree, we have that $h_{G_i}(v_{i+j}) = h_{G_{i+j}}(v_{i+j})$
for $j \le 3$. Consequently, $v'_l = h_K(v)$, which finishes the proof.
\qed
\end {proof}

Combining Theorems~\ref{niceTT} and~\ref{wniceTT} we obtain 
a corollary.

\begin {corollary}   \label{corniceTT}
An undirected nice graph is left $M$-homotens for every ring~$M$.
A (directed or undirected) weakly nice graph 
is left $M$-homotens for every ring~$M \ne \zet_2^k$.
\end {corollary}

Extending our conditions that guarantee that a graph 
is $M$-homotens, we present the following lemma, which will
be used in Section~\ref{sec:applications}. Note that the
assumption about spanning subgraphs is needed. 

\begin {lemma}   \label{suphomotens}
Suppose $H$ contains a connected spanning $M$-homotens
graph. Then $H$ is $M$-homotens.
\end {lemma}

\begin {proof}
Let $f : H \TT[M] K$, let $G$ be the connected spanning 
$M$-homotens subgraph of~$H$. Restriction of~$f$ to~$E(G)$
is $TT_M$, hence $f(e) = g^\sharp(e)$ for each $e \in E(G)$
and some (anti)homomorphism~$g$. Let $e = uv \in E(H) \setminus E(G)$.
We have to prove $f(e) = (g(u),g(v))$. Let $P$ be a path from~$u$
to~$v$ in~$G$. By treating the closed walk $P \cup \{ uv \}$
as $W$ in the end of the proof of Theorem~\ref{wniceTT}, we conclude
the proof.
\qed
\end {proof}

\subsection{Applications} \label{sec:applications}

In this section we provide several applications of nice graphs
(that is of Theorem~\ref{wniceTT} and Corollary~\ref{corniceTT}).
Particularly, we prove that `almost all' graphs are left $M$-homotens
for every ring~$M$ and construct an embedding of category~$\calG_{hom}$ into
$\calG_{TT_M}$. This result was proved (for $M=\zet_2$)
in~\cite{NS-TT1} by an ad-hoc construction. Here we follow a more systematic
approach---we employ a modification of an edge-based replacement
operation (see \cite{HN}). 
As a warm-up we prove an easy, but perhaps surprising result.

\begin {corollary}   \label{nicesupgraph}
For every graph~$G$ there is a graph~$G'$ containing $G$ as an
induced subgraph such that for every ring~$M$ every 
$TT_M$~mapping from~$G'$ to arbitrary graph is induced by a homomorphism
(i.e., $G'$~is $M$-homotens).
\end {corollary}

\begin {proof}
We take as $G'$ the (complete) join of~$G$ and $K_5$; that is, we let
$V(G') = V(G) \cup \{v_1,v_2,\dots, v_5\}$, and
$E(G') = E(G) \cup \{\mbox{all edges containing some $v_i$}\}$.
By Theorem~\ref{wniceTT} it is enough to show that $G'$ is nice. Every 
copy of $K_t$ ($t < 5$) in~$G'$ can be extended to~$K_5$ by adding
some vertices $v_i$. One can also show routinely that any two copies 
of~$K_4$ in~$G'$ are `$K_4$-connected'---condition~4 in 
Definition~\ref{def:wnice}.
\qed
\end {proof}

The following theorem was our main motivation for
introducing (weakly) nice graphs. Note that `a.a.s.' means, as usual, 
`asymptotically almost surely', that is `with probability 
tending to 1'.

\begin {theorem}   \label{homotens}
Let $M$ be a ring.
\begin {enumerate}
  \item Complete graph $K_k$ is $M$-homotens for $k \ge 5$
     (and for $k\ge 4$ if $M \ne (\zet_2)^t$).
  \item The random graph $G(n,1/2)$ is $M$-homotens a.a.s.
  \item The random $k$-partite graph is $M$-homotens a.a.s. for $k \ge 5$
     (and for $k\ge 4$ if $M \ne (\zet_2)^t$). Explicitly, 
     $$
        \lim_{n\to \infty} \Pr[ \hbox{$G = G(n,1/2)$ is $M$-homotens} 
         \mid \hbox{$G$ is $k$-partite} ] = 1 \,.
     $$
  \item The random $K_k$-free graph is $M$-homotens a.a.s. for $k \ge 6$
     (and for $k\ge 5$ if $M \ne (\zet_2)^t$).
\end {enumerate}
If $M \ne \zet_2$, then in each of the statement, any orientation
of the considered graph is $M$-homotens, too.
\end {theorem}

\begin {proof}
As $K_t$ is nice (weakly nice for $t = 4$), 1 follows by
Corollary~\ref{corniceTT}.
In~\cite{NS-TT1} we proved that the random graph is a.a.s. nice, 
so again, Corollary~\ref{corniceTT} implies 2. 
By~\cite{KPR}, a random $K_k$-free graph is a.a.s. $(k-1)$-partite, 
hence 3 implies 4. The proof of 3 is similar to the proof that the
random graph is a.a.s. nice, we sketch it for convenience.

Let $A_1$, \dots, $A_k$ be the parts of the random $k$-partite graph.
By standard arguments, all $A_i$'s are a.a.s. approximately
of the same size, in particular all are non-empty. It is a routine to 
verify parts 1, 2, and (in case $k \ge 5$) 3 of Definition~\ref{def:wnice}.
For part 4, let $V(K) = \{v_1, \dots, v_4\}$, $V(K') = \{v_9, \dots, v_{12}\}$. 
We pick $i_1$, \dots, $i_4$ so that $v_t \not\in A_{i_k}$, except
possibly if $t = k$ or $t=k+8$. 
We attempt to pick $v_5 \in A_{i_1}$, \dots, $v_8 \in A_{i_4}$ to
satisfy the condition~4. The probability that a particular 4-tuple
fails is at most $\bigl(1-2^{-18}\bigr)^{n/2k}$.
Hence, the probability that some copies $K$, $K'$ of~$K_4$
are `bad' is at most $n^8 c^n$ (for some $c < 1$).
\qed
\end {proof}

We proceed by another application of Corollary~\ref{corniceTT}
--- we show that the structure of $TT_M$~mappings is at least
as rich as that of homomorphisms.

\begin {theorem}   \label{niceembedding}
There is a mapping $F$ that assigns graphs to graphs, such that 
for any ring~$M$ and for any graphs $G$,~$H$ 
(we stress that we consider loopless graphs only) holds
$$
   G \leh H  \iff   F(G) \leTT_M F(H) \,.
$$
Moreover $F$ can be extended to a 1-1 correspondence for mappings
between graphs: 
if $f:G \to H$ is a homomorphism, then $F(f):F(G)\to F(H)$
is a $TT_M$~mapping and any $TT_M$~mapping
between $F(G)$ and~$F(H)$ is equal to~$F(f)$ for some 
homomorphism~$f : G \hom H$.
(In category-theory terms, $F$ is an embedding of the category of
all graphs and their homomorphisms into the category of
all graphs and all $TT_M$-mappings between them.)
\end {theorem}

\begin {figure}
\centerline{
  \includegraphics{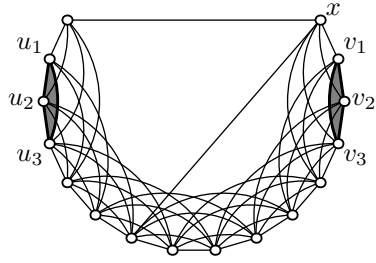}
}
  \caption{The graph~$I$ used in triangle-based replacement
  (proof of Theorem~\ref{niceembedding}).}
  \label{fig:indicator}
\end {figure}

\begin {proof}
We will use a modification of edge-based replacement 
(see \cite{HN}). Let $I$ be the graph in Figure~\ref{fig:indicator}
with arbitrary (but fixed) orientation.
To construct $F(G)$, we will replace each of the 
vertices of~$G$ by a triangle and each of the edges of~$G$
by a copy of~$I$, gluing different copies on triangles. 
More precisely, let $U = V(G)\times{}\{0,1,2\}$, for every 
edge $e \in E(G)$ let $I_e$ be a separate copy of~$I$. 
If $e = (u,v)$ then we identify vertex $u_i$ ($i \in \{0,1,2\}$)
with $(u,i)$ in~$U$, and vertex $v_i$ with $(v,i)$ in~$U$.
Let $F(G)$ be the resulting graph; we write shortly $F(G) = G*I$.
If $f: V(G) \to V(H)$ is a homomorphism then 
we define $F(f):E(F(G)) \to E(F(H))$ as follows: let
$e=(u,v)$ be an edge of~$G$ and $a$ an edge of $E(I_e)$. 
Let $e'$ be the image of~$e$ under~$f$. In the isomorphism 
between $I_e$ and $I_e'$ the edge $a$ gets mapped to some~$a'$.
We put $F(f)(a) = a'$. It is easily seen that $F(f)$ is a $TT_\zet$
(thus $TT_M$) mapping that is induced by a homomorphism, 
we let $\phi(f)$ denote this homomorphism.
Now, we turn to the more difficult step of proving that 
every $TT_M$ mapping from~$G$ to~$H$ is~$F(f)$ for
some $f:G \hom H$. We will need several auxiliary claims.

\textbf{Claim 1.} $I$ is critically $6$-chromatic. 

Take any $K_5$ in~$I$, color in by 5 colors. There is a unique
way how to extend it, which fails, so $\chi(I) \ge 6$. 
Clearly 6 colors suffice. Moreover, if we delete any vertex of~$I$
then it is possible to color the remaining vertices consecutively
$1, 2, 3, 4, 5, 1, 2, \dots, 5$.

\textbf{Claim 2.} $I$ is rigid.

That is, the only homomorphism $f: I \to I$ is the identity. 
By Claim~1, $f$ cannot map $I$ to its subgraph, hence $f$ is
an automorphism. There is a unique vertex~$x$ of degree~9, so
$f$ fixes it. There is a unique hamiltonian cycle $x=x_1, \dots, x_?$
such that $x_i x_j$ is an edge whenever $|i-j| \le 4$, therefore
this cycle has to be fixed by~$f$ too. This leaves two possibilities, 
but only one of them maps the edge $x_1 x_?$ properly.

\textbf{Claim 3.} $I$ is $K_5$-connected.

That is, for every two vertices $a$, $b$ of~$I$ there is a path
$a=a_1, a_2, \dots,a_k = b$ such that $a_i a_j$ is an edge whenever 
$|i-j| \le 4$. 

\textbf{Claim 4.} Whenever $H$ is a graph and $g:I \hom H*I$
a homomorphism, there is an edge $e \in E(H)$ such that $g$
is an isomorphism between $I$ and~$I_e$. 

If $g$ maps all vertices of $I$ to one of the $I_e$'s, then
we are done by Claim~2. If not, let $a$, $b$ be vertices of~$I$
such that $g(a)$ is a vertex of~$I_e$ (for some edge $e=uv\in E(H)$)
and~$g(b)$ is not. Choose a path $a=a_1, a_2, \dots,a_k = b$ as in Claim~3.
Let $a_i$ be the last vertex on this path that is a vertex of~$I_e$.
Not all three vertices $a_{i-1}$, $a_{i-2}$, $a_{i-3}$ can be 
in the `connecting triangle' $\{v\}\times{}\{0,1,2\}$, on the other
hand each of them is connected to $a_{i+1}$, a contradiction.

\textbf{Claim 5.} For every graph~$H$ the graph $H*I$ is nice. 

This is an easy consequence of Lemma~\ref{suphomotens}.

To finish the proof, let $h: F(G)\TT F(H)$ be a $TT_M$ mapping. As
graph $G*I$ is nice, it is $M$-homotens by Corollary~\ref{corniceTT}.
Therefore $h$ is induced by a homomorphism, say $g: F(G) \hom F(H)$.
By Claim~4, $g$ maps an~$I_e$ to an~$I_{e'}$, therefore there is a
homomorphism $f:V(G) \to V(H)$ such that $g = \phi(f)$ and $h = F(f)$,
as claimed. 
\qed
\end {proof}

\subsection{A necessary condition}
\label{sec:necessary}

In this section we present a necessary condition for a graph
to be $\zet$-homotens.%
\footnote{As any $TT_\zet$ mapping is 
  $TT_M$ for every ring~$M$ (Lemma~\ref{groupinfluence}),
  each $M$-homotens graph is also $\zet$-homotens. Therefore, 
  the presented condition is necessary for a graph to be
  $M$-homotens, too.
  To illustrate that $\zet_2$-homotens is indeed stronger condition
  than $\zet$-homotens, we note that no 4-chromatic graph is
  $\zet_2$-homotens---it admits a $TT_2$~mapping to $K_3$.
}
As mentioned earlier, odd circuits are the simplest examples of 
graphs that are not $\zet$-homotens.
Similarly, no graph with a vertex of degree 2 is $\zet$-homotens, 
except of a triangle.
This way of thinking can be further strengthened and generalized, 
yielding Theorem~\ref{chromcon}.
To state our result in a compact way, we introduce a definition
from~\cite{GNN}. We say that a graph~$G$ is \emph{chromatically 
$k$-connected} if for every $U \subseteq V(G)$ such that $G - U$ is
disconnected the induced graph $G[U]$ has chromatic number at least~$k$.
Equivalently~\cite{GNN}, $G$ is chromatically $k$-connected, iff
every homomorphic image of~$G$ is $k$-connected.

\begin {theorem}   \label{chromcon}
Let $M$ be a ring. If a graph is connected and $M$-homotens then it
is chromatically $3$-connected. 
\end {theorem}

\begin {proof}
Suppose $G$ is a counterexample to the theorem. 
Hence, vertices of~$G$ can be partitioned into 
sets $A$, $B$, $U$, $L$, such that $A \cup B$
separates $U$ from~$L$; that is there is no edge
from~$U$ to~$L$, moreover $A$, $B$ are independent sets.
We may suppose $A \cup B$ is a minimal set that separates
$U$ from~$L$.
We are going to prove that $G$~is not $\zet$-homotens, therefore 
by Lemma~\ref{groupinfluence} not $M$-homotens as well.

We identify all vertices of~$A$ to a single vertex~$a$, and
all vertices of~$B$ to a vertex~$b$. Let $F$ be the resulting graph,
and $f:G \to F$ be the identifying homomorphism.
We define a $TT_\zet$~mapping $g$ from~$F$ as follows.
For $u \in U$ we map edge $(u,a)$ (if it exists) to~$(b,u)$, 
$(a,u)$ to~$(u,b)$, $(u,b)$ to~$(a,u)$, and $(b,u)$ to~$(u,a)$.
For $u, v \in U$ we map edge $(u,v)$ (if it exists) to~$(v,u)$.
Every other edge is mapped to itself. We let $F'$~denote
the resulting graph (it has the same set of vertices as~$F$).
It is straightforward to use Lemma~\ref{altdef} to 
verify that $g$~is indeed $TT_\zet$.

Hence $gf^\sharp$ is a $TT_\zet$~mapping; we need to show that it is
not induced.  At least one of~$A$, $B$ is non-empty. Suppose it is~$A$
and pick $x \in A$. As $A \cup B \setminus \{x\}$ is not a separating
set ($A \cup B$ is a minimal one), there are vertices $u \in U$ and 
$l \in L$ that are adjacent to~$x$, without loss of generality $(x,u)$,
$(x,l)$ are edges of~$G$.  By definition of~$g$ we have
$gf^\sharp((x,l)) = (x,l)$ and $gf^\sharp((x,u)) = (u,y)$. Therefore
$gf^\sharp$ maps two adjacent edges to two nonadjacent edges, hence it
is not induced.  
\qed
\end {proof}

The following corollary deduces a simpler necessary condition, 
though a weaker one: We can prove that the graph of icosahedron
is not $\zet$-homotens by using Theorem~\ref{chromcon} (the 
neighborhood of an edge is a $C_6$), but not using
Corollary~\ref{triangle}.

\begin {corollary}   \label{triangle}
Let $G$ be a connected graph with at least four vertices. Suppose
the neighbourhood of some $v \in V(G)$ induces a bipartite graph.
Then $G$ is not $M$-homotens for any ring~$M$.

Consequently, every vertex of a homotens graph 
is incident with an odd wheel (in particular with a triangle), 
except if it is contained in a component of size at most three.
\end {corollary}

\begin {proof}
Let $A$, $B$ be the color-classes of neighborhood of~$v$. 
If there is a vertex nonadjacent to~$v$, then we can use 
Theorem~\ref{chromcon}. So suppose $v$ is connected to every
vertex of~$G$. Then every other vertex has a bipartite neighborhood.
The only case that stops us from using Theorem~\ref{chromcon} is
when $|A|,|B|\le 1$, that is when $G$~has at most three vertices.
\qed
\end {proof}

A somewhat surprising consequence of Corollary~\ref{triangle} is that
no triangle-free graph is homotens. This immediately
answers a question of~\cite{NS-TT1}. It also implies, that a connected cubic
graph is $M$-homotens only if it is a $K_4$ and $M$ is not a power
of~$\zet_2$. More generally, we have the following result
(compare Theorem~\ref{homotens}).

\begin {corollary}   \label{norandom}
Let $r\ge 3$ be integer, $M$ ring. 
The probability that a random $r$-regular graph is $M$-homotens 
is bounded by a constant less than 1, if size of the graph
is large enough.
\end {corollary}

\begin {proof}
It is known \cite{Wormald-survey} that the probability
that random $r$-regular graph is triangle-free tends
to a nonzero limit, hence we can apply Theorem~\ref{triangle}.
\qed
\end {proof}

Corollary~\ref{triangle} also indicates that complete graphs
involved in the definition of nice graphs are necessary, at least to
some extent. 
However, the condition of Corollary~\ref{triangle} (or
Theorem~\ref{chromcon}) is far from
being sufficient: for example the graph from Proposition~\ref{cnhomotens}
is chromatically $k$-connected and not $\zet$-homotens.
In particular, we do not know whether there are $K_4$-free homotens graphs.
By~\cite{KPR}, a random $K_4$-free graph is a.a.s. $3$-partite, hence
not chromatically 3-connected, hence by Theorem~\ref{chromcon}
not $\zet$-homotens. Still, it is possible that $K_4$-free
$\zet$-homotens graphs exist, promising candidates are Kneser graphs
$K(4n-1,n)$, which are chromatically 3-connected for large~$n$
\cite{GNN}.

\begin {question}   \label{q:homotens}
Is the Kneser graph $K(4n-1,n)$ $\zet$-homotens, if $n$ is large enough?
\end {question}

\section{Right homotens graphs} \label{sec:rh}

In this section we complement Section~\ref{sec:lh} by study
of graphs which, when used as target graphs, make existence
of $TT$~mappings and of homomorphisms coincide.
Recall that a graph~$H$ is called right $M$-homotens if the existence
of a $TT_M$~mapping from an arbitrary graph to~$H$ implies the existence
of a homomorphism. Right homotens graphs (in comparison with 
left homotens ones) provide more structure; in this section we 
characterize them by means of special Cayley graphs
and state a question aiming to find a better characterization.

\subsection{Free Cayley graphs}
\label{sec:delta}

Free Cayley graphs were introduced by Naserasr and~Tardif~\cite{NaserasrTardif} 
(see also thesis of Lei Chu \cite{LeiChu})
in order to study chromatic number of Cayley graphs.
They will serve us as a tool to study $TT$~mappings, 
in particular we will use them to study right homotens graphs
and to prove density in Section~\ref{sec:density}.

Let $M$ be a ring, let $H$ be a graph.
For a vertex $v \in V(H)$ we let $e_v : V(H) \to M$ be the indicator function, 
that is $e_v(u) = 1$ if $v=u$ and $e_v(u) = 0$ otherwise.
We define graph\footnote{
  More precisely, we define $\Delta_M(H)$ to be a directed graph.
  However, if $\symH$~is a symmetric orientation of an undirected
  graph~$H$, then $\Delta_M(\symH)$ is a symmetric orientation of some
  undirected graph~$H'$, we may let $\Delta_M(H) = H'$. 
  The whole Section~\ref{sec:delta} may be modified for undirected graphs
  by similar changes.}
$\Delta_M(H)$ 
with vertices~$M^{V(H)}$, where $(f,g)$ is an edge iff 
$g - f = e_v - e_u$ for some edge $(u,v) \in E(H)$.
We can see that $\Delta_M(H)$ is a Cayley graph, 
it is called the \emph{free Cayley graph} of~$H$.
We begin our study of free Cayley graphs with a simple observation and
with a useful lemma, which is due to Naserasr and Tardif 
(for a proof, see~\cite{LeiChu}).

\begin {proposition}   \label{trivemb}
Graph $\Delta_M(H)$ contains $H$ as an induced subgraph.
\end {proposition}

\begin {proof}
Take functions $\{e_v \mid v \in V(H)\} \subseteq V(\Delta_M(H))$.
\qed
\end {proof}

\begin {lemma}   \label{fcg}
Let $M$ be a ring, $H$ a Cayley graph on~$M^k$ (for some integer~$k$)
and $G$ an arbitrary graph. Then any homomorphism $G \hom H$
can be (uniquely) extended to a mapping $\Delta_M(G) \to H$
that is both graph and ring homomorphism. 
\end {lemma}

The following easy lemma appears in~\cite{DNR} (although without
explicit mention of graphs~$\Delta_M$).

\begin {lemma} \label{delta}
$G \TT[M] H$ is equivalent with $G \hom \Delta_M(H)$.
\end {lemma}

Note that Lemma~\ref{cube} is a special case of Lemma~\ref{delta}, 
as graphs $\Delta(G)$ defined in Section~\ref{sec:ex}
are isomorphic to $\Delta_{\sss \zet_2}(G)$. 
Lemmas~\ref{fcg} and~\ref{delta} have as immediate corollary an
embedding result that nicely complements Theorem~\ref{niceembedding}.
In contrary with Theorem~\ref{niceembedding} though,
our embedding is not functorial, it is just embedding of quasiorder
$(\calG, \leTT_M)$ in $(\calG, \leh)$.

\begin {corollary}   \label{deltaembedding}
$G \TT[M] H$ is equivalent with $\Delta_M(G) \hom \Delta_M(H)$.
\end {corollary}

\begin {proof}
If $G \TT[M] H$ then by Lemma~\ref{delta} we have
$G \hom \Delta_M(H)$ and by Lemma~\ref{fcg} the result
follows. For the other implication, by Proposition~\ref{trivemb}
graph $G$ maps homomorphically to $\Delta_M(H)$, and 
application of Lemma~\ref{delta} yields $G \TT[M] H$.
\qed
\end {proof}

We remark that Corollary~\ref{deltaembedding} provides an embedding of
category of $TT_M$~mappings to category of Cayley graphs with mappings
that are both ring and graphs homomorphisms.

\subsection{Right homotens graphs}

We start with two simple observations concerning right homotens
graphs. The first one is a characterization of 
right homotens graphs by means of~$\Delta_M$. 
It does not, however, give an efficient method (polynomial algorithm)
to verify if a given graph is right homotens, neither a good
understanding of right homotens graphs. Hence, we will seek better
characterizations (compare with Corollary~\ref{rhequiv} and
Question~\ref{rhquestion}).

\begin {proposition}   \label{rhdelta}
A graph $H$ is right $M$-homotens if and only if $\Delta_M(H) \hom H$.
\end {proposition}

\begin {proof}
For the `only if' part it is enough to observe that $\Delta_M(H) \TT[M] H$
for every graph $H$: clearly $\Delta_M(H) \hom \Delta_M(H)$ and
we use Lemma~\ref{delta}.
For the other direction, if $G \TT[M] H$ then 
by Lemma~\ref{delta} we have $G \hom \Delta_M(H)$ and by composition 
(Lemma~\ref{compose}) we have $G \hom H$. 
\qed
\end {proof}

\begin {lemma}   \label{rhimpl}
Assume $H \hom H'$ and $H' \TT[M] H$. If $H$ is right $M$-homotens 
then $H'$ is right $M$-homotens as well.
\end {lemma}

\begin {proof}
If $H$ is right $M$-homotens, then $\Delta_M(H) \hom H$.  
By Corollary~\ref{deltaembedding} from $H' \TT[M] H$ we deduce that
$\Delta_M(H') \hom \Delta_M(H)$. By composition, 
$$
   \Delta_M(H') \hom \Delta_M(H) \hom H \hom H' \,,
$$
hence $H'$ is right $M$-homotens. 
\qed
\end {proof}

\begin {corollary}   \label{rhequiv}
Let $H$, $H'$ be homomorphically equivalent graphs (that is
$H \hom H'$ and~$H' \hom H$). Then $H$ is right $M$-homotens
if and only if $H'$ is right $M$-homotens.
\end {corollary}

Note that $TT_M$-equivalence is not sufficient in Corollary~\ref{rhequiv}:
each graph $H$ is $TT_M$-equivalent with $\Delta_M(H)$
and the latter is always a right $M$-homotens graph (for each~$M$), as
we will see from the next proposition. 
Also note that the analogy of Corollary~\ref{rhequiv} does not
hold for left homotens graphs.

Next, we consider a class of right $M$-homotens graphs that is 
central to this topic. 
We will say that $H$ is an \emph{$M$-graph} if it is a Cayley graph on
some power of~$M$ ($\zet_2$-graphs are also called cube-like
graphs; they have been introduced by Lov\'asz~\cite{Harary} as an
example of graphs, for which every eigenvalue is an integer).

\begin {proposition}   \label{rhcbl}
Any $M$-graph is right $M$-homotens.
\end {proposition}

\begin {proof}
Let $H$ be an $M$-graph. As $H \hom H$, by Lemma~\ref{fcg} we conclude
that $\Delta_M(H) \hom H$. 
\qed
\end {proof}

In analogy with the chromatic number we define the $TT_M$~number 
$\chiTT[M](G)$ to be the minimum~$n$ for which
there is a graph~$H$ with $n$ vertices such that
$G \TT[M] H$. As any homomorphism induces a $TT_M$ mapping, 
we see that $\chiTT[M](G) \le \chi(G)$ for every graph~$G$.
Continuing our project of finding similarities between $TT_M$~mappings and
homomorphisms, we prove that for finite~$M$ the $TT_M$~number cannot be much
smaller than the chromatic number.

\begin {corollary}   \label{rhcomplete}
Let $G$ be arbitrary graph.
If $M$ is a finite ring of characteristic $p$ then
$\chi(G) / \chiTT[M](G) < p$.

Moreover, $\chi(G) / \chiTT[\zet](G) < 2$.
\end {corollary}

\begin {proof}
First we prove that $\chi(G) < m \cdot \chiTT[M](G)$
for any finite ring $M$ of size $m$.
To this end, consider a Cayley graph on $M^k$ with the generating
set $M^k \setminus \{\vecz\}$---that is a complete graph~$K_{m^k}$ 
with every edge in both orientations. This is
an $M$-graph, hence by Proposition~\ref{rhcbl} it is right $M$-homotens. 

Now, choose $k$ so that $m^{k-1} < \chiTT[M](G) \le m^k$. It follows
that
$G \TT[M] K_{m^k}$, and as $K_{m^k}$ is right $M$-homotens,
$G \hom K_{m^k}$. Therefore,
$\chi(G) \le m^k < m \cdot \chiTT[M](G)$.

Next, if $p$ is the characteristic of~$M$, this means that $M$
contains~$\zet_p$ as a subring. This by Lemma~\ref{groupinfluence}
implies that any $TT_M$ mapping is $TT_{\zet_p}$, thus
$\chiTT[M](G) \ge \chiTT[\zet_p](G)$, and the result follows.
For the second part we use Lemma~\ref{groupinfluence} again to 
infer that any $TT_\zet$ mapping is $TT_{\zet_2}$.
\qed
\end {proof}

How good is the bound given by Corollary~\ref{rhcomplete}
is an interesting and difficult question. Even in the simplest
case $M = \zet_2$ this is widely open; perhaps surprisingly this is
related with the quest for optimal error correcting codes. For
details, see~\cite{NS-TT1,RS-thesis}
Another corollary of Proposition~\ref{rhcbl} is a characterization
of right homotens graphs.

\begin {corollary}   \label{rhchar}
A graph is right $M$-homotens if and only if it is homomorphically
equivalent to an $M$-graph.
\end {corollary}

\begin {proof}
The `if' part follows from Lemma~\ref{rhequiv} and
Statement~\ref{rhcbl}. For the `only if' part, notice that $\Delta_M(H)$
is a $M$-graph, $H \subseteq \Delta_M(H)$, and if $H$ is right
$M$-homotens then $\Delta_M(H) \hom H$.
\qed
\end {proof}

Corollary~\ref{rhchar} is not very satisfactory, as it does not
provide any useful algorithm to verify if a given graph is right
homotens. Indeed, it is more a characterization of graphs
that are hom-equivalent to some $M$-graph, than the other way around:
Suppose we are to test if a given graph is hom-equivalent to some (arbitrarily
large) $M$-graph. It is not obvious if there is a finite process that
decides this; however Corollary~\ref{rhchar} reduces this task to 
decide if $\Delta_M(H) \hom H$. 
The latter condition is easily checked by an obvious brute-force algorithm.

We hope that a more helpful characterization of right homotens graphs will 
result from considering the core of a given graph. As a core of a graph~$H$ 
is hom-equivalent with~$H$, it is right homotens if and only if $H$~is.
Therefore, we attempt to characterize right homotens cores, leading to an
easy proposition and an adventurous question. We note that one part of the
proof of the proposition is basically the folklore fact that the core of a
vertex-transitive graph is vertex-transitive, while the other part is a
generalization of an argument used by~\cite{HKSS} to prove that $K_n$
is right $\zet_2$-homotens if and only if $n$ is a power of~2.
However, we include the proof for the sake of completeness.

\begin {proposition}   \label{rhcores}
Let $H$ be a right $M$-homotens graph that is a core. Then 
\begin {itemize}
  \item $|V(H)|$ is a power of~$|M|$, and
  \item $H$ is vertex transitive. If $M=\zet_2$,
     then for every two vertices of~$H$, there is an automorphism
     exchanging them.
\end {itemize}
\end {proposition}

\begin {proof}
For a function $g \in M^{V(H)}$ we let $H_g$ denote the subgraph 
of~$\Delta_M(H)$ induced by the vertex set $\{g + e_v; v \in V(H)\}$.
Observe that each $H_g$ is isomorphic with~$H$. Let $f: \Delta_M(H) \to H$
be a homomorphism and for each $u\in V(H)$, define 
$V_u = \{ v\in V(\Delta_M(H)); f(v) = u \}$. 
Now $f$ restricted to $H_g$ is a homomorphism
from~$H_g$ to~$H$. As $H$~is a core, every homomorphism from~$H$
to~$H$ is a bijection. Consequently, for every~$g$
the graph $H_g$ contains precisely one vertex from each $V_u$. 
By considering all graphs $H_g$ we see that all sets $V_u$ are of
the same size $|M|^{|V(H)|}/|V(H)|$. Therefore, $|V(H)|$ is a power
of~$|M|$.

For the second part let $u$, $v$ be distinct vertices of~$H$. 
We know $\Delta_M(H) \hom H$. As $H \simeq H_\vecz$
($\vecz$ being the identical zero), we have 
a homomorphism $f: \Delta_M(H) \hom H_\vecz$. As $H$ is a core, we
know that $f$ restricted to $H_\vecz$ is an automorphism
of~$H_\vecz$. By composition with the inverse automorphism, we may 
suppose that $f$ restricted to $H_\vecz$ is an identity.
Next, consider the isomorphism 
$\phi: \Delta_M(H) \hom \Delta_M(H)$ given by $g \mapsto g + e_v - e_u$.
A composed mapping $f \circ \phi$ is a homomorphism $H_\vecz \hom
H_\vecz$ (therefore an automorphism) that maps $u$ to~$v$. 
Moreover, if $M = \zet_2$ then $f \circ \phi$ maps $v$ to~$u$
as well.
\qed
\end {proof}

The previous proposition suggests that a stronger result might be true, 
and that this may be a way to a characterization of right homotens
graphs. In particular, we ask the following.

\begin {question}   \label{rhquestion}
\begin {enumerate}
  \item Suppose $H$ is a right $M$-homotens graph and a core. 
    Is $H$ an $M$-graph?
  \item Is the core of each $M$-graph an $M$-graph?
\end {enumerate}
\end {question}

We note that even the (perhaps easier to understand) case $M = \zet_2$
is open. 
But one can see easily that 1 and~2 in Question~\ref{rhquestion} are equivalent:
If $H$ is a right $M$-homotens core, then $H$ is the core of the 
$M$-graph~$\Delta(H)$; hence 2 implies 1.
Conversely, let $K$ be an $M$-graph and $H$ its core.
By Proposition~\ref{rhcbl}, $K$ is right $M$-homotens,
therefore by Corollary~\ref{rhequiv} $H$~is right $M$-homotens.
If 1 is true, then $H$ is an $M$-graph, as claimed.

\section{Density} \label{sec:density}

In this section we compare homomorphisms and tension-continuous 
mappings from a different perspective: we prove that partial orders
defined by existence of a homomorphism (a $TT_M$~mapping respectively) 
share an important property, namely the \emph{density}. 
To recall, we say that a partial order $<$ is dense, if for every $A$, $B$ 
satisfying $A < B$ there is an element~$C$ for which $A < C < B$.

It is known \cite{HN,NT} that the homomorphism order (with all 
hom-equivalence classes of finite graphs as elements
and with the relation $\lh$) is dense, if we do not consider
graphs without edges. The parallel result for the order
defined by $TT_M$~mappings is given by the following theorem. 
In fact we prove a stronger property (proved in~\cite{HN}
for homomorphisms) that every finite antichain in a given
interval can be extended; density is the special case $t=0$.

\begin {theorem}   \label{AntiExtInt}
Let $M$ be a ring, let $t \ge 0$ be an integer.
Let $G$, $H$ be graphs such that $G \lTT_M H$ and $E(G) \ne \emptyset$.
Let $G_1$, $G_2$, \dots, $G_t$ be pairwise
incomparable (in $\lTT_M$) graphs satisfying $G\lTT_M G_i \lTT_M H$
for
every $i$. Then there is a graph~$K$ such that 
\begin {enumerate} 
  \item $G \lTT_M K \lTT_M H$,
  \item $K$ and~$G_i$ are $TT_M$-incomparable for every $i = 1, \dots, t$.
\end {enumerate}
If in addition $G \leh H$ then we have even $G \lh K \lh H$.
If we consider undirected graphs, then we get undirected graph~$K$. 
\end {theorem}

This theorem was proved in a previous paper~\cite{NS-TT1} by the authors,
here we present a much shorter proof. The key of the proof is
the use of graphs~$\Delta_M(G)$ for a new proof of Lemma~\ref{SIL}.
From this, Theorem~\ref{AntiExtInt} follows directly.

\begin {proof}[Theorem~\ref{AntiExtInt}---sketch]
We use the next lemma for graphs $G$, $G_1$, \dots, $G_t$.
and we let $G'$ be the graph, that this lemma ensures.
Put $K = G + G'$. For details, see \cite{NS-TT1}.
\qed
\end {proof}

\begin {lemma}[Sparse incomparability lemma for $TT_M$]    \label{SIL}
Let $M$ be an abelian group (not necessarily a finitely generated one), 
let $l$, $t\ge 1$ be integers.
Let $G_1$, $G_2$,~\dots, $G_t$, $H$ be (finite directed non-empty\footnote{
  that is with non-empty edge set}) 
graphs such that $H \nTT[M] G_i$ for every~$i$.
Then there is a graph~$G$ such that
\begin {enumerate}
  \item $g(G)>l$ (that is $G$~contains no circuit of size at most~$l$), 
  \item $G \lh H$,
  \item $G \nTT[M] G_i$ for every $i = 1, \dots, t$.
\end {enumerate}
(For undirected graphs we get undirected graph~$G$.)
\end {lemma}

In the proof we will use a variant of Sparse incomparability lemma 
for homomorphisms in the following form (it has been proved for
undirected graphs in~\cite{NZ}, the version we present here follows
by the same proof).

\begin {lemma}[Sparse incomparability lemma for homomorphisms] \label{homoSIL}
Let $l$, $t\ge 1$ be integers,
let $H$, $G_1$, \dots, $G_t$ be (finite directed non-empty) graphs such that
$H \nhom G_i$ for every $i$. Let $c$ be an
integer. Then there is a (directed) graph~$G$ such that 
\begin {itemize}
  \item $g(G)>l$ (that is $G$~contains no circuit of size at most~$l$),
  \item $G \lh H$, and
  \item $G \nhom G_i$ for every~$i$.
\end {itemize}
(For undirected graphs we get undirected graph~$G$.)
\end {lemma}

Before we start the proof, we summarize necessary results about influence
of ring~$M$ on the existence of $TT_M$ mappings. The following 
summarizes results that appear as Theorem~4.4 in \cite{DNR}, and as
Lemma~14 and~17 in~\cite{NS-TT1}.

\begin {lemma}   \label{groupinfluence}
Let $G$, $H$ be graphs, $f: E(G) \to E(H)$ any mapping.
\begin {enumerate}
  \item If $f$ is $TT_\zet$ then it is $TT_M$ for any group~$M$.
  \item Let $M$ be a subring of $N$. If $f$ is $TT_N$ then it is $TT_M$.
  \item Let $G$, $H$ be \emph{finite} graphs.
      Then $G \leTT_n H$ holds either for finitely many~$n$ or 
      for every $n$. In the latter case $G \leTT_{\sss \zet} H$ holds.
\end {enumerate}
\end {lemma}

\begin {proof}[Lemma~\ref{SIL}]
First, suppose that $M$ is a finite ring;
by Lemma~\ref{delta} we know that $H \nhom \Delta_M(G_i)$ for every~$i$.
Therefore, we may use Lemma~\ref{homoSIL} to obtain $G'$ of girth greater
than~$l$ such that $G' \leh H$ and $G' \not\leh \Delta_M(G_i)$.
Consequently $G' \nTT_M G_i$ for every~$i$.  

Next, let $M$ be an infinite, finitely generated group, that is a ring.
Then $M \simeq \zet^\alpha \times{}\prod_{i=1}^k \zet_{n_i}^{\beta_i}$, 
for some integers $k$, $n_i$, $\beta_i$, $\alpha$. As $M$ is
infinite, we have $\alpha > 0$, therefore $M \ge \zet$.
By Lemma~\ref{groupinfluence} we conclude that for any mapping it is 
equivalent to be $TT_M$ and to be $TT_\zet$, hence we may suppose
$M = \zet$. By Lemma~\ref{groupinfluence}, there is only finitely
many integers $n$ for which holds $H \TT[n] G_i$ for some~$i$ or 
$H \TT[n] \edge$. Pick some $n$ for which neither of this holds.
By the previous paragraph for ring $\zet_n$ we find a graph~$G'$ 
such that $G' \nTT[n] G_i$ for every $i = 1, \dots, t$.
It follows from Lemma~\ref{groupinfluence} that also $G' \nTT[M] G_i$.

Finally, let $M$ be a general abelian group.
For each mapping $f: E(H) \to X$ (where
$X \in \{G_1, \dots, G_t\}$) there is 
an $M$-tension~$\phi_X$ on~$X$ which certifies that $f$ is not
a $TT_M$ mapping. Let 
  $A = \big\{ \phi_X(e) \mid 
                 e \in E(X), X \in \{G_1, \dots, G_t\} \big\}$ 
be the set of all elements of~$M$ that are used for these certificates.
Let $M'$ be the subgroup of~$M$ generated by~$A$; by the choice
of~$A$ we have $H \nTT[M'] G_i$.
By the previous paragraph there is a graph~$G'$ that meets 
conditions 1, 2, and $G' \nTT[M'] G_i$ for every~$i$.
Consequently, $G' \nTT[M] G_i$ for every~$i$, which concludes
the proof.
\qed
\end {proof}

Let us add a remark that partially explains the way we conducted the
above density proof.
Standard proofs of density of the homomorphism order rely on the fact, 
that the category of graphs and homomorphisms has products. 
We prove next, that this is not true for $TT_M$~mappings; 
therefore another approach is needed.
In~\cite{NS-TT1} we developed a new structural Ramsey-type theorem 
to overcome the non-existence of products;
here we used the construction $\Delta_M$ for much shorter proof. 

\begin {proposition}   \label{noproducts}
Category $\calG_{TT_M}$ of (directed or undirected) graphs 
and $TT_M$ mappings does not have products for any ring~$M$.
\end {proposition}

\begin {proof}
We will formulate the proof for the undirected version, although
for the directed version the same proof goes through.
We show that there is no product $C_3 \times{}C_3$.  Suppose, to
the contrary, that $P$~is the product $C_3 \times{}C_3$. 
Let $\pi_1$, $\pi_2: P \TT C_3$ be the projections, 
let $E(C_3) = \{e_1, e_2, e_3\}$. 

We look first at mappings $f_i: \edge \to C_3$ sending the
only edge of $\edge$ to $e_i$. If we consider mapping $f_i$ to the first 
copy of~$C_3$ and $f_j$ to the second one, by definition of
the product there is exactly one edge~$e \in E(P)$ such that 
$\pi_1(e) = e_i$ and $\pi_2(e) = e_j$. We let $e_{i,j}$ denote this~$e$.
So, $E(P)$ consists of nine edges $e_{i,j}$, for $1 \le i, j \le 3$.

As $\pi_1$, $\pi_2$ are $TT$ mappings, by Lemma~\ref{altdef}
there are no loops in~$P$. There are no parallel edges either: 
suppose $e$, $f$ are parallel edges in~$P$. Then without loss
of generality $\pi_1(e) \ne \pi_1(f)$, hence we get a contradiction
by Lemma~\ref{altdef}.

Finally, for a $\rho \in S_3$ let $f_\rho: C_3 \to C_3$
send $e_i$ to $e_{\rho(i)}$. Using the definition of product for mapping
$f_{\id}$ and $f_{\rho}$, Lemma~\ref{altdef}, and the fact that
there are no parallel edges in~$P$ we find that 
$E_\rho = \{e_{1,\rho(1)}, e_{2,\rho(2)}, e_{3,\rho(3)}\}$ are edges
of a cycle.  Considering $\rho = \id$ and $\rho=(1,3,2)$ we find that
part of~$P$ looks as in the Figure~\ref{fig:catproof}
(in the directed case, the orientation may be arbitrary, 
if $M = \zet_2^k$).

\begin{figure}[h]
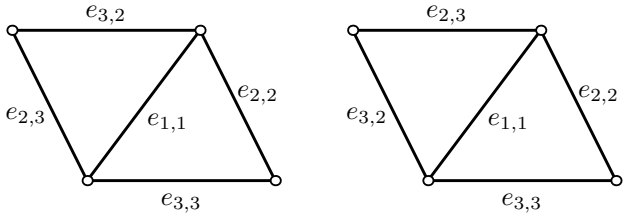

\centerline{
  \hfil
  \includegraphics{noproducts.0}
  \hfil
  \includegraphics{noproducts.1}
  \hfil
}
\caption{Proof of Proposition~\ref{noproducts}.}
\label{fig:catproof}
\end{figure}

Consider the first case. As $E_\rho$ is a cycle for $\rho=(2,3,1)$, 
the edges~$e_{1,2}$ and~$e_{2,3}$ are adjacent. By taking
$\rho=(2,1,3)$, we find that $e_{1,2}$ and~$e_{2,3}$ are adjacent.
As there are no parallel edges in~$P$, we have $e_{1,2} = xy$~or
$e_{1,2} = yx$. Hence, $e_{1,2}$, $e_{2,3}$, $e_{2,2}$ forms a cycle. 
As $\pi_1$ is $TT$ mapping, we obtain a contradiction by
Lemma~\ref{altdef}. In the second case we proceed in the same way with
edge~$e_{2,1}$, we prove that it is adjacent with $e_{3,2}$
and~$e_{3,3}$ and yield a contradiction with $\pi_2$ being a $TT$~mapping.
\qed
\end {proof}

\section {Remarks}

\subsection{Broader context (Jaeger's project)}

Tension-continuous mappings were defined in~\cite{DNR,NS-TT1} 
in a broader context of three related types of mappings:
$FF$ (lifts flows to flows), $FT$ (lifts tensions to flows), 
and $TF$ (lifts flows to tensions). 
In \cite{DNR,RS-thesis} these mappings are studied in more detail, 
in particular their connections to several classical 
conjectures (Cycle Double Cover conjecture, Tutte's 5-flow conjecture, and
Berge-Fulkerson matching conjecture) are explained. 

The universality and density of $TT$ mappings shows that the
Jaeger's project of characterizing ``atoms'' of a partial order
defined by flow-continuous mappings has no dual analogue (for $TT$
mappings). It follows from Theorem~\ref{AntiExtInt} that
each of the quasiorders $\leTT_M$ is everywhere dense for the class
of directed graphs. Graphs $\orK_2$ and the loop graph are
the minimal and the maximal element of these orders. Particularly, 
there cannot be any atom (the contrary is conjectured for the
flow-continuous order in~\cite{DNR,Jaeger}). This is also in sharp
contrast with the homomorphism order of oriented graphs where the
homomorphism order~$\leh$ contains many gaps of a complicated
structure. (These gaps are characterized by~\cite{NT}.)
Another consequence of Theorem~\ref{AntiExtInt} is that each of the
orders $\leTT_M$ contains an infinite antichain, a property which
is presently open for $M$-flow-continuous mappings for every~$M$, 
in particular for cycle-continuous mappings; see~\cite{DNR}.

\subsection{TT-perfect graphs}

For every graph $G$, its chromatic number $\chi(G)$ is at least as big
as the size of its largest clique, $\omega(G)$. Recall, that a graph~$G$ is
called \emph{perfect} if $\chi(G') = \omega(G')$ holds for every
induced subgraph~$G'$ of~$G$. A graph is called \emph{Berge} if
for no odd $l \ge 5$ does $G$ contain $C_l$ or $\barC_l$ as an
induced subgraph. It is easy to see that being perfect implies
being Berge; the so-called Strong Perfect Graph Conjecture (due to
Claude Berge) claims that the opposite is true, too.
Perfect graphs have been a topic of intensive research 
that recently lead to a proof \cite{CRSTperfect} of the 
Strong Perfect Graph Conjecture.

As a humble parallel to this development we define a graph~$G$ to
be \emph{$TT$-perfect}\footnote{%
  more precisely, $TT_2$-perfect, 
  but we will not consider $M \ne \zet_2$ in this section
}
if for every induced subgraph~$G'$ of~$G$
we have $\chiTT[2](G') \le \omega(G')$
(definition of $\chiTT[2](G')$ appears 
before Corollary~\ref{rhcomplete}).
Equivalently, $G$ is $TT$-perfect if each of its induced
subgraphs~$G'$ admits a $TT_2$ mapping to its maximal clique. 

Note that we cannot ask for $\chiTT(G') = \omega(G')$ 
since $K_4 \TT K_3$, and therefore
$\chiTT(K_4) = 3$, while $\omega_{\sss TT}(K_4) = 4$.

As any homomorphism induces a $TT$ mapping (see Lemma~\ref{homo}), 
$\chiTT(G') \le \chi(G')$ holds for every graph~$G'$.
Consequently, 
every perfect graph is $TT$ perfect. The converse, however, is 
false. For example, let $G = \barC_7$. Graph $G$ itself is 
not perfect.
On the other hand $\chiTT(G) = 3$ and every induced
subgraph of $G$ is Berge, hence perfect, hence $TT$-perfect.
Let us study $TT$-perfect graphs in a similar manner
as Strong Perfect Graph Theorem does for perfect graphs. 
To this end, we define a graph~$G$ to be \emph{critical} 
if $G$ is not $TT$-perfect, but each induced subgraph of~$G$ is.
We start our approach by a technical lemma.

\begin {lemma}   \label{TT_pefect_circuits}
Let $l \ge 3$ be odd. Cycle $C_l$ is not $TT$-perfect. 
Graph $\barC_l$ is $TT$-perfect if and only if $l = 7$.
\end {lemma}

\begin {proof}
Clearly $\chiTT(C_l) = 3 > \omega(C_l)$. Graph $\barC_7$
was discussed above, $\barC_5$ is isomorphic to~$C_5$. 
As $\chi(\barC_9) = 5$ and as $K_4$ is right $\zet_2$-homotens, 
being a $\zet_2$-graph, we have
$\chiTT(\barC_9) = 5 >  \omega(\barC_9)$. 
It is easy to verify that graphs~$\barC_l$ for $l \ge 13$
are nice. Thus they are homotens and not $TT$-perfect, since
they are not perfect. The only remaining case is the 
graph~$\barC_{11}$. This is not nice, on the other hand,
every edge is contained it a $K_5$ and all $K_5$'s are 
`connected'---there is a chain of all 11 copies of $K_5$
such that neighboring copies intersect in a~$K_4$.
It follows that $\barC_{11}$ is homotens, in particular
$\barC_{11} \not\TT K_5$.
\qed
\end {proof}

\begin {figure}
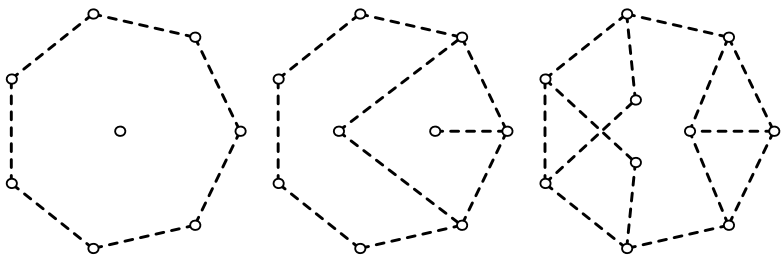
   
\centerline{
  \hfil
  \includegraphics{critical.0}
  \hfil
  \includegraphics{critical.1}
  \hfil
  \includegraphics{critical.2}
  \hfil
}
\caption{Several critical graphs that are not cycles neither
complements of cycles. The dashed lines denote precisely 
the {\bf non-edges} of the graph.
}
\label{fig:critical}
\end {figure}

\begin {corollary}   \label{critical}
For every odd $l>3$ graph $C_l$ is critical; if 
$l \ne 7$ then $\barC_l$ is critical, too.
Moreover graphs $G_1$, $G_2$, and $G_3$ in Figure~\ref{fig:critical}
are critical.
\end {corollary}

\begin {proof}
We sketch the proof of $G_1$ being critical. We have
$\chi(G_1) = 1 + \chi(\barC_7) = 5$, therefore Corollary~\ref{rhcomplete}
implies $\chi_{TT}(G_1) = 5 > \omega(G_1)$ and $G_1$ is not $TT$-perfect.
Let $G'$ be an induced subgraph of~$G_1$. If $G' = \barC_7$
then $G'$ is $TT$-perfect; otherwise, it is a routine to verify that
$G'$~is Berge, consequently perfect and $TT$-perfect.
\qed
\end {proof}

We do not know how many other critical graphs there are, not even if there
is an infinite number of them. 

\section* {Acknowledgements}

The authors would like to thank Ji\v r\'\i\ Matou\v sek for a
stimulating question that lead to the notion of right homotens graphs
and to L\'aszlo Lov\'asz for pointing us to~\cite{Shih-thesis},
which was an inspiration for Proposition~\ref{cnhomotens}.


\bibliographystyle{rs-amsplain}
\bibliography{pariz}

\providecommand{\bysame}{\leavevmode\hbox to3em{\hrulefill}\thinspace}
\providecommand{\MR}{\relax\ifhmode\unskip\space\fi MR }
\providecommand{\MRhref}[2]{%
  \href{http://www.ams.org/mathscinet-getitem?mr=#1}{#2}
}
\providecommand{\href}[2]{#2}
\begin{thebibliography}{10}

\bibitem{LeiChu}
Lei Chu, \emph{Colouring {C}ayley graphs}, Master's thesis, University of
  Waterloo, 2004.

\bibitem{CRSTperfect}
Maria Chudnovsky, Neil Robertson, Paul Seymour, and Robin Thomas, \emph{The
  strong perfect graph theorem}, Annals of Mathematics, to appear.

\bibitem{DNR}
Matt DeVos, Jaroslav~Ne\v set\v ril, and Andr\'e Raspaud, \emph{On flow and
  tension-continuous maps}, KAM-DIMATIA Series \textbf{567} (2002).

\bibitem{Diestel}
Reinhard Diestel, \emph{Graph theory}, Graduate Texts in Mathematics, vol. 173,
  Springer-Verlag, New York, 2000.

\bibitem{GNN}
Chris~D. Godsil, Richard~J. Nowakowski, and Jaroslav Ne{\v{s}}et{\v{r}}il,
  \emph{The chromatic connectivity of graphs}, Graphs Combin. \textbf{4}
  (1988), no.~3, 229--233.

\bibitem{HKSS}
Ge{\v n}a Hahn, Jan Kratochv{\'\i{}}l, Jozef {\v S}ir{\'a}{\v n}, and Dominique
  Sotteau, \emph{On the injective chromatic number of graphs}, Discrete Math.
  \textbf{256} (2002), no.~1-2, 179--192.

\bibitem{Harary}
Frank Harary, \emph{Four difficult unsolved problems in graph theory}, Recent
  advances in graph theory (Proc. Second Czechoslovak Sympos., Prague, 1974),
  Academia, Prague, 1975, pp.~249--256.

\bibitem{HN}
Pavol Hell and Jaroslav~Ne\v set\v ril, \emph{Graphs and homomorphisms}, Oxford
  Lecture Series in Mathematics and Its Applications, Oxford University Press,
  2004.

\bibitem{Jaeger}
Fran{\c{c}}ois Jaeger, \emph{On graphic-minimal spaces}, Ann. Discrete Math.
  \textbf{8} (1980), 123--126, Combinatorics 79 (Proc. Colloq., Univ.
  Montr\'eal, Montreal, Que., 1979), Part I.

\bibitem{Kelmans}
Alexander~K. Kelmans, \emph{On edge bijections of graphs}, Tech. Report 93-41,
  DIMACS, 1993.

\bibitem{KPR}
Phokion~G. Kolaitis, Hans-J\"urgen Pr{\"o}mel, and Bruce~L. Rothschild,
  \emph{{$K\sb {l+1}$}-free graphs: asymptotic structure and a {$0$}-{$1$}
  law}, Trans. Amer. Math. Soc. \textbf{303} (1987), no.~2, 637--671.

\bibitem{KR}
V{\'a}clav Koubek and Vojt{\v e}ch R{\"o}dl, \emph{On the minimum order of
  graphs with given semigroup}, J. Combin. Theory Ser. B \textbf{36} (1984),
  no.~2, 135--155.

\bibitem{LMT}
Nathan Linial, Roy Meshulam, and Michael Tarsi, \emph{Matroidal bijections
  between graphs}, J. Combin. Theory Ser. B \textbf{45} (1988), no.~1, 31--44.

\bibitem{NaserasrTardif}
Reza Naserasr and Claude Tardif, \emph{Chromatic numbers of {C}ayley graphs on
  $\zet_2^n$}, manuscript.

\bibitem{Nes-derivative}
Jaroslav Ne{\v{s}}et{\v{r}}il, \emph{Homomorphisms of derivative graphs},
  Discrete Math. \textbf{1} (1971), no.~3, 257--268.

\bibitem{NT}
Jaroslav Ne{\v{s}}et{\v{r}}il and Claude Tardif, \emph{Duality theorems for
  finite structures (characterising gaps and good characterisations)}, J.
  Combin. Theory Ser. B \textbf{80} (2000), no.~1, 80--97.

\bibitem{NZ}
Jaroslav Ne{\v{s}}et{\v{r}}il and Xuding Zhu, \emph{On sparse graphs with given
  colorings and homomorphisms}, J. Combin. Theory Ser. B \textbf{90} (2004),
  no.~1, 161--172, Dedicated to Adrian Bondy and U. S. R. Murty.

\bibitem{RS-thesis}
Robert {\v S}\'amal, \emph{On {XY} mappings}, Ph.D. thesis, Charles University,
  2006.

\bibitem{NS-TT1}
Jaroslav~Ne\v set\v ril and Robert {\v S}\'amal, \emph{Tension-continuous
  maps---their structure and applications}, submitted, arXiv:math.CO/0503360.

\bibitem{Shih-thesis}
Ching-Hsien Shih, \emph{On graphic subspaces of graphic spaces}, Ph.D. thesis,
  The Ohio State University, 1982.

\bibitem{Whit1}
Hassler Whitney, \emph{{Congruent graphs and the connectivity of graphs}}, Am.
  J. Math. \textbf{54} (1932), 150--168.

\bibitem{Wormald-survey}
Nicholas~C. Wormald, \emph{Models of random regular graphs}, Surveys in
  combinatorics, London Math. Soc. Lecture Note Ser., vol. 267, Cambridge Univ.
  Press, Cambridge, 1999, pp.~239--298.

\end{thebibliography}

\end{document}